\documentclass[11pt]{amsart}
\usepackage{amsmath}
\usepackage{amssymb}
\usepackage{epsfig}
\usepackage{a4wide}

\newcommand{\B}{\ensuremath{\mathbb{B}}}
\newcommand{\N}{\ensuremath{\mathbb{N}}}
\newcommand{\Z}{\ensuremath{\mathbb{Z}}}
\newcommand{\R}{\ensuremath{\mathbb{R}}}
\newcommand{\T}{\ensuremath{{\mathcal T}}}
\newcommand{\bs}[1]{\boldsymbol{#1}}
\newcommand{\V}{\ensuremath{\boldsymbol V}}

\DeclareMathOperator{\inn}{int}
\DeclareMathOperator{\cl}{cl}
\DeclareMathOperator{\supp}{supp}
\DeclareMathOperator{\diam}{diam} 

\newtheorem{thm}{Theorem}[section] 
\newtheorem{lemma}[thm]{Lemma}
\newtheorem{defi}[thm]{Definition}
\newtheorem{corollary}[thm]{Corollary} 

\allowdisplaybreaks[3]

\begin{document} 

\title{Computing Modular Coincidences for substitution tilings and
  point sets} 

\author{D.~Frettl\"oh}
\author{B.~Sing}
\address{Fakult\"at f\"ur Mathematik, Universit\"at
  Bielefeld, Postfach 100131, 33501 Bielefeld,  Germany}
\email{frettloe@math.uni-bielefeld.de, sing@math.uni-bielefeld.de}
\urladdr{http://www.math.uni-bielefeld.de/baake}

\begin{abstract} Computing modular coincidences can show whether
  a given substitution system, which is supported on a point lattice
  in  $\R^d$, consists of model sets or not. We prove the
  computatibility of this problem  and determine an upper bound 
  for the number of iterations needed. The main tool is a simple
  algorithm for computing modular coincidences, which is essentially a 
  generalization of Dekking coincidence to more than one dimension,
  and the proof of equivalence of this generalized Dekking coincidence
  and  modular coincidence. As a consequence, we also obtain some
  conditions for the existence of modular
  coincidences. In a separate section, and throughout the article, a
  number of examples are given. 

\end{abstract} 

\maketitle

\section{Introduction}\label{sec:intro}
Substitution tilings, or substitution point sets, are examples of objects
which show high order at both 
short and long range, even if they are not periodic. This manifests
itself in the diffraction spectrum of such a tiling. Like
crystals or periodic tilings, some nonperiodic substitution tilings show
a diffraction pattern consisting of bright spots --- the ``Bragg peaks''
--- only. Such tilings are called \emph{pure point} diffractive. In
general, it is hard to decide whether a given substitution tiling is
pure point diffractive. For a special class of tilings the question can
be answered for any given example in principle: If the set of
difference vectors between tiles spans a 
discrete point lattice, then there is a condition which is
equivalent to pure point diffractivity (\cite{LMS03}, see also Theorem
\ref{thm:modcoinc} below).  This class covers many
well-known examples, such as the chair tiling and the sphinx tiling
(being pure point, see, e.g., \cite{Sol97} \cite{BMS98}, \cite{Rob99},
\cite{LM01}), or the Thue-Morse sequence and the table tiling 
(being not pure point, see, e.g., \cite{Mah26}, \cite{Kak72},
\cite{Rob99}). In fact, the condition of~\cite{LMS03} may be difficult to
check in some cases. Here, we devise a simple algorithm and 
show the equivalence between the output of this algorithm and the
condition mentioned. As a consequence, we obtain results about the
computability of the problem and some conditions for the existence
(resp.~nonexistence) of modular coincidences. 

The paper is organized as follows: Section \ref{sec:basics} contains
notations and basic terms used in this text, the definition of lattice
substitution systems and related facts and definitions.  Section
\ref{sec:modsets} is dedicated to the main result of \cite{LMS03} (Theorem
\ref{thm:modcoinc}), together with the definitions of model sets and
modular coincidence. Since we used the notation of \cite{LMS03} wherever
it is possible, readers who are familiar with \cite{LMS03} may skip
Sections \ref{sec:basics} and \ref{sec:modsets}. In Section
\ref{sec:main}, we give the construction of the coincidence graph,
which is a generalization of Dekking coincidence \cite{Dek78}, 
and establish the equivalence between this coincidence and modular
coincidence (Theorem \ref{thm:aequiv}). From this result, we obtain 
results about the computability of the problem (Theorems
\ref{thm:upperbd}, \ref{thm:vmldv}), as well as  
some conditions (Corollaries \ref{cor:psising}, \ref{cor:paircoinc},
\ref{cor:bij}) for the existence (resp.~nonexistence) of 
modular coincidences. Section \ref{sec:complex} is a brief remark
about the best upper bound for the costs of computing modular
coincidences. Section \ref{sec:examples} is a collection of
examples which illustrate terms and statements of preceding sections
and show the usefulness of our results.   
Finally, three appendices show the connections between this text and
Dekking coincidence, the overlap algorithm in \cite{Sol97} and automatic
sequences, respectively.

\section{Basic Definitions, Substitution Systems}\label{sec:basics}
Let us fix some notation. By $\N$ we denote the set of all positive
integers. The sum of two sets is to be understood point-wise, i.e.,
$M+N=\{ n+m \mathbin| m\in M, n\in N\}$. The closure of a set $M$ is
denoted by $\cl(M)$, the interior of $M$ by $\inn(M)$. The cardinality of $A$ 
is denoted by $\#A$ (where $\# (\varnothing) = 0$). By $\B^d$ we denote the
closed unit ball $\{ x \mathbin| \| x \| \le 1 \}$ in $\R^d$. 
A \emph{tile}  is a nonempty compact set $T \subset \R^d$ with the property
that  $\cl(\inn(T))=T$. A \emph{tiling} of $\R^d$ is a collection of tiles 
$\{T^{}_n\}^{}_{n\in\N}$ which covers $\R^d$ and contains no overlapping
tiles, i.e., $\inn(T^{}_k) \cap \inn(T^{}_n) = \varnothing$ for $k \ne n$.   
A tiling $\T$ is called \emph{nonperiodic}, if the only solution of 
$\T+x=\T$  is $x=0$. If the tiles of a tiling belong to
finitely many translation classes $[T^{}_1],[T^{}_2], \ldots, [T^{}_m]$, 
then the class representatives  $T^{}_1,T^{}_2,\ldots,T^{}_m$ are called
\emph{prototiles}. If necessary, one may consider two congruent
prototiles as different by labelling them with their type. 
Given a set of prototiles $T_i$, one can describe tilings --- or parts
of tilings --- by a set 
\begin{equation} \label{eq:tplusx}
 \{ T^{}_1 + v^{(1)}_1, T^{}_1 + v^{(1)}_2, \ldots, T^{}_2 + v^{(2)}_1 ,
T^{}_2+ v^{(2)}_2, \ldots , T^{}_m+v^{(m)}_1, \ldots \} \qquad
(v^{(j)}_i \in \R^d),
\end{equation}  
or shortly by $\{T^{}_1 + V^{}_{1}, T^{}_2 + V^{}_{2} \ldots, T^{}_m +
V^{}_{m}\}$ --- where $V^{}_{i}$ contains all
$v^{(i)}_k$ of \eqref{eq:tplusx}) --- or even shorter by a
\emph{position vector} $\V=(V^{}_{1},V^{}_{2}, \ldots, V^{}_{m})^T$.

Nonperiodic tilings of long range order can easily be generated by a
substitution. The idea is, basically, to give a set of prototiles 
and a rule how to enlarge them, and then to dissect them into tiles 
congruent to the original prototiles (cf.~Example
\ref{subsec:chairtable}). Formally, this can be done by a 
\emph{matrix function system}, resp.~by a \emph{substitution system}. 

\begin{defi} \label{def:mfs}
A \emph{matrix function system (MFS)} is an $m \times m$--matrix
$\Phi=(\Phi^{}_{ij})_{1 \le i,j \le m}$, where each $\Phi^{}_{ij}$ is a
finite (possibly empty) set of mappings $x \mapsto Q x + a^{}_{ijk}$, 
where $Q$ is an expansive linear map, $1 \le i,j  \le m$, $1 \le k \le \#
(\Phi^{}_{ij})$, $a^{}_{ijk} \in \R^d$. \\
An MFS $\Phi$ is called \emph{primitive} if its
\emph{substitution matrix} $S(\Phi)=(\#(\Phi^{}_{ij}))^{}_{1\le i,j\le
  m}$ (an $m\times m$-matrix with non-negative integer entries) is
primitive,  i.e., $S^{k}_{}$ has only strictly positive entries for
some $k\in\N$.   
\end{defi}

Note that we can define $\Phi^{k}_{}$ (for $k\in\N$) inductively by the
composition of maps, i.e., \linebreak
$\Phi^{k+1}_{}=\Phi^{k}_{}\circ\Phi=(\Phi^{k}_{}\circ\Phi)^{}_{1\le
i,j\le m}=(\bigcup_{\ell=1}^{m} (\Phi^{k})_{i\ell}\circ(\Phi^{})_{\ell
j})^{}_{1\le i,j\le m}$ where
\begin{equation*}
(\Phi^{k})_{i\ell}\circ\Phi^{}_{\ell j} = \left\{ \begin{array}{ll}
\varnothing & \textnormal{if }(\Phi^{k})_{i\ell}=\varnothing
\textnormal{ or }\Phi^{}_{\ell j}=\varnothing, \\
\{g\circ f\mathbin| g\in(\Phi^{k})_{i\ell},\,f\in\Phi^{}_{\ell j}\} &
\textnormal{otherwise.} \\
\end{array}\right.
\end{equation*}
Clearly, we have $S(\Phi^{k}_{})\le S(\Phi)^{k}_{}$. Note also that
the linear part of the affine mappings in $\Phi^{k}_{}$ is given by
$x\mapsto Q^{k}_{}x$. Now we can apply the MFS $\Phi$ to a position
vector $\V$. I.e., if $\V=(V_1,\ldots,V_m)$ is a position vector as
above, then $\Phi \V = (W_1, \ldots, W_m)$, 
$W_i= \bigcup_{j=1}^m \Phi^{}_{ij}(V^{}_j) = \bigcup_{j=1}^m
\bigcup_{f\in\Phi^{}_{ij}} f(V^{}_j)$.
 
\begin{defi}\label{def:sub_sys}
A \emph{substitution system} is a pair $(\V,\Phi)$, where
$\V=(V^{}_{1},V^{}_{2}, \ldots, V^{}_{m})^T$ is a fixed point of the 
matrix function system $\Phi$ (i.e., $\Phi\V=\V$), such that all
$V^{}_{i}$ are mutually disjoint and the unions in
\begin{equation*}
V^{}_i = \bigcup_{j=1}^m \Phi^{}_{ij}(V^{}_j) = \bigcup_{j=1}^m
\bigcup_{f\in\Phi^{}_{ij}} f(V^{}_j)
\end{equation*} 
are pairwise disjoint for all $1\le i\le m$.\\
If $\operatorname{supp}(\V)= \bigcup_{i=1}^m V^{}_i$ is a lattice,
then $(\V,\Phi)$ is called \emph{lattice substitution  system (LSS)
  with $m$ components}. 
\end{defi}

Now the relation between substitution systems and substitution tilings 
becomes clear: A set of tiles is represented by a vector
$\V$ whose entries $V_i$ are point sets (where each point encodes the
position of one tile). The substitution applied to this set of
tiles yields a larger set of tiles, which is represented by the vector
$\Phi \V$.

By the action of a tile substitution, each tile $T$ is enlarged by
$Q$, then dissected into some tiles of type $[T^{}_1],\ldots,
[T^{}_m]$. If $T \in [T^{}_i]$, then $QT$ is dissected into
$\#(\Phi^{}_{1i})$ tiles congruent to $T^{}_1$, $\#(\Phi^{}_{2i})$
tiles congruent to $T^{}_2$ and so on. The translational part of the
affine maps in $\Phi$ determines the exact location of these tiles
inside the enlarged tile. In this way, $(\T,\Phi)$ can be considered a
substitution system. Wherever we will emphasize the fact that
$(\T,\Phi)$ describes a tiling, we will call $\T$ a 
\emph{substitution tiling} and $\Phi$ a tile substitution. 
For examples of the relation between tile substitutions and
substitutions on point sets, we refer to Example 
\ref{subsec:chairtable}. The next lemma reads essentially: The
translate of an LSS $(\V,\Phi)$ by a $t \in \supp(\V)$ is again an
LSS.  

\begin{lemma} \label{lem:lsstrans}
Let $(\V,\Phi)$ be an LSS, $t \in \supp(\V)$,
and let $(\V-t,\Phi+Qt-t)$ be defined by 
$\V-t = (V_1 - t, \ldots, V_m - t)$ and $\Phi+Qt-t =
(\Phi_{ij}+Qt-t)_{ij}, \; \Phi_{ij}+Qt-t=  \{ \varphi+Qt-t \mathbin| 
\varphi \in \Phi_{ij} \}$. Then $(\V-t,\Phi+Qt-t)$ is also an LSS.  
\end{lemma}
\begin{proof}
A simple calculation
(cf.~Def.~\ref{def:sub_sys}) yields: $(\Phi+Qt-t)(\V-t)=\V-t$ if and
only if $\Phi(\V)=\V$. Since $t \in \supp(\V)$, it follows that
$\supp(\V-t)$ is a lattice, thus all conditions for lattice substitution
systems are satisfied. 
\end{proof}

If $(\V,\Phi)$ is an LSS, then 
the set $\V$ in Definition \ref{def:sub_sys} consists of point sets
$V_i$. Rather than considering $\V$ as a subset of $\R^{dm}$, we may
think of $\V$ as a point set in $\R^d$, where each point has got
one \emph{color} (out of $m$ colors altogether). Such a single point 
of $\V$ with color $j$, located at $x \in \R^d$,  is of the form 
\begin{equation*}
\{(\varnothing, \ldots, \varnothing, \underbrace{\{ x \}
}_{j\textnormal{th entry}}  , \varnothing, \ldots,\varnothing ) \}, 
\end{equation*}
which is a somewhat inconvenient notation. Therefore, we will avoid to
consider such single points as much as possible. If such a single
colored point located at $x$ of type $j$ is considered, we use the
notation '$\V \cap \{x\}$, $x$ is of type $j$' (or of color $j$).

In this paper we consider only primitive LSS
$(\V,\Phi)$. This implies that each $V^{}_i$ is a \emph{Delone set},
and we call $\V$ a \emph{multi-component Delone set}
(see~\cite[Section 2]{LMS02}, where the name \emph{Delone multiset} is
used; however, the name ``multiset'' in this context is misleading).

\begin{defi}
A set $V \subset \R^d$ is a \emph{Delone set} (or \emph{Delaunay set}) if
it satisfies the following two conditions:
\begin{enumerate}
\item $V$ is \emph{uniformly discrete}: There is a value $r$ such that each
  open ball of radius $r$ contains at most one point of $V$.
\item $V$ is \emph{relatively dense}: There is a value $R$ such that each
  closed ball of radius $R$ contains at least one point of $V$.
\end{enumerate}
\end{defi}
Note that the maximal such radius $r$ is the \emph{packing radius} for
$V$, while the minimal value of $R$ is the \emph{covering radius} of $V$.

We should mention that every primitive substitution tiling gives rise
to a primitive substitution system in a quite natural way, see
\eqref{eq:tplusx} and comments there. In the other direction, this is less
obvious. In the case of LSS, one simple way is the following: Assign
to each point $x$ of type $i$ in a multi--component Delone set a tile
$x+F$ of type $i$, where $F$ is a fundamental domain of the lattice.   
Then it is trivial to read off the substitution for the tiling from
the MFS for the multi--component Delone set. Note that we do not
require the tile--substitution to be self--similar in the sense of
\cite{Sol97}. For details about the subtleties of the relation between
point--substitutions and tile--substitutions, see \cite[Section
3.2]{LMS03}. The property which describes this relation  
is given in the following definition. It works both for tilings and
for multi--component Delone sets. In order to express it, 
we denote by $C_r(x,A)$ the finite portion of the structure $A$
in a ball of radius $r$ around $x$, the \emph{$r$--cluster} at $x$.
If $A$ is a tiling, then $C_r(x,A):= \{ T \in A \mathbin| T
\cap (r \B^d +x) \ne \varnothing \}$. If $A=(A_1,\ldots,A_m)$ is a
multi-component Delone set, then $C_r(x,A):= ( A_1 \cap (r \B^d +x),
\ldots , A_m \cap (r \B^d +x) )$.
 
\begin{defi} \label{def:mld}
Let $A,B$ be two multi--component Delone sets, or two tilings, or one
multi--component Delone set and one tiling. 
$B$ is \emph{locally derivable (with radius $r$)} from $A$, if, for all
$x,y$ in $\R^d$, one has: 
\begin{equation*}
C_r(x,A)=C_r(y,A) + (x-y) \quad \Longrightarrow \quad C_0(x,B)=C_0(y,B) + (x-y).
\end{equation*}   
If $B$ is locally derivable from $A$ and $A$ is locally derivable from
$B$, then $A$ and $B$ are \emph{mutually locally derivable (MLD)},
see~\cite{BSJ91}.  
\end{defi}

For examples of a tiling and its corresponding substitution system
being MLD see Example \ref{subsec:chairtable}. Another property
which is frequently used is the following.

\begin{defi}
A tiling $A$, or a multi--component Delone multiset $A$, is called
\emph{repetitive}, if for every $r>0$ there is $R>0$ such that each
$R$--cluster in $A$ contains a translate of each $r$--cluster in $A$.    
\end{defi}

\section{Model Sets}\label{sec:modsets}

The objects we consider in this paper are primitive lattice
substitution systems $(\V,\Phi)$
and tilings which are MLD to a primitive LSS.
To exclude certain pathological cases, we assume further that all
$r$--clusters in $\V$ are contained in some $\Phi^k({\bs v})$, where
$k \in \N$ and ${\bs v}$ is a single point of some type $i$. In
particular, this ensures that $\V$ is repetitive. Moreover, 
the dynamical system associated to $\V$ is always minimal, see the
remark at the end of this section or \cite{Sol97} for details. 

A central question is whether such a 
substitution system consists of model sets or not; and whether the
diffraction spectrum of such a substitution pattern is pure point or
not. If $\V$ is fully periodic (i.e., there are $d$ linearly
independent translation vectors such that $\V$ is invariant under
these translations), this is trivial. The question becomes interesting
if $\V$ is not fully periodic. 
A theorem by \textsc{Schlottmann}~\cite[Theorem 4.5]{Sch00} states
that all regular model sets are pure point diffractive.
We first recall what a model set is, see~\cite{Moo00}.

\begin{defi} \label{def:modset}
A \emph{cut and project scheme} consists of a collection of spaces and
mappings as follows:
\begin{equation*}
\begin{array}{ccccl}
\R^d & \stackrel{\pi}{\longleftarrow} & \R^d\times H &
\stackrel{\pi^{}_{\textnormal{int}}}{\longrightarrow} & H \\
\cup & & \cup & & \cup \textnormal{ dense}\\
L & \stackrel{1-1}{\longleftarrow} & \tilde{L} & \longrightarrow & L^{\star} \\
\end{array}
\end{equation*}
where $H$ is a locally compact Abelian group, $\tilde{L}\subset
\R^d\times H$ is a lattice, i.e., a discrete subgroup for which the
quotient group $(\R^d\times H)/\tilde{L}$ is compact, and $\pi$,
$\pi^{}_{\textnormal{int}}$ denote the canonical projections, where
$\pi^{}_{\textnormal{int}}(\tilde{L})=L^{\star}$ is dense in $H$ and
the restriction $\pi|^{}_{\tilde{L}}$ is injective (and therefore
  bijective as map from $\tilde{L}$ to $L=\pi(\tilde{L})$).
We call $\R^d$ the \emph{direct space} (or \emph{physical space}) and
$H$ the \emph{internal space}.
\\
A \emph{model set} (or \emph{cut and project set}) in $\R^d$ is a
subset of $\R^d$ which (up to translation) is of the form
$\varLambda(\Omega)=\{\pi(x)\mathbin| x\in\tilde{L},\,
\pi^{}_{\textnormal{int}}(x)\in\Omega\}$ for some cut and project
scheme, where $\Omega\subset H$ has non-empty interior and compact
closure. We call $\Omega$ the \emph{window} of the model set.\\
A model set $\varLambda(\Omega)$ is \emph{regular} if the boundary
$\partial\Omega=\cl(\Omega)\setminus\inn(\Omega)$ is of Haar measure $0$.   
\end{defi}
Note that one can define the so-called \emph{star map} $x_{}^{\star}$
of $x\in L$ by $x^{\star} =
\pi^{}_{\textnormal{int}}(\pi|^{-1}_{\tilde{L}}(x))$, justifying the
notation $L^{\star}$ above.

To an LSS, one can naturally associate a cut
and project scheme (see Definition~\ref{def:modcoinc} below for the
definition of $L'$). The internal space is then given by 
$H=(L/L')\times \lim_{\leftarrow k} L'/Q^k L'$, where $\lim_{\leftarrow
  k} L'/Q^k L'$ denotes the $Q$-adic completion of $L'$ (a profinite
group), see~\cite[Section 5.1]{LMS03}, and Examples
\ref{subsec:perseq} and \ref{subsec:abab}. Furthermore, the 
lattice $\tilde{L}$ also arises naturally as diagonal embedding of
$L$ in $\R^d\times H$, i.e., $\tilde{L}=\{(t,t)\mathbin| t\in L\}$. So
only the existence and shape of a window is in question.
For this, a central role is played by the notion of a modular coincidence.

\begin{defi}[cf.~\cite{LM01,LMS03}]\label{def:modcoinc}
Let $(\V,\Phi)$ be an LSS and
$L=\operatorname{supp}(\V)$. Let $L':= 
L_1 + L_2 + \cdots +L_m$, where $L_i:= \langle V_i - V_i \rangle_{\Z}$. \\
For $a \in L$, let 
\begin{equation*}
\begin{split}
\Phi^{}_{ij}[a] 
& := \{ \varphi \in  \Phi^{}_{ij} \mathbin| \varphi(y) \equiv a \bmod
QL',\textnormal{ where } V^{}_j \subset y+L' \} \\ 
& \; = \{ \varphi \in \Phi^{}_{ij} \mathbin|
  \varphi(V^{}_j)\subset a+QL'\}. 
\end{split}
\end{equation*}
Furthermore, let 
\begin{equation*}
\Phi[a] := \bigcup_{1 \le i,j \le m} \Phi^{}_{ij}[a]. 
\end{equation*}
$(\V,\Phi)$ admits a \emph{modular coincidence relative to $QL'$},
if $\Phi[a]$ is contained entirely in one row of $\Phi$ for some
$a\in L$, i.e., if $(a+QL')\subset V^{}_i$ for some $1\le i\le m$. 
\end{defi}

For an illustration of this concept, see Example
\ref{subsec:kol24} with an explicit calculation of a modular
coincidence in a particular case. (At this point, consider just the
first part of the example, before the first occurrence of $\Psi$).  

\textsc{Remark:} Since $Q(L) \subseteq L$, one also has
$Q(L') \subseteq L'$ (see~\cite[Section
5.1]{LMS03}). Furthermore, $L/L'$ is finite, i.e., $L'$ partitions $L$
in finitely many cosets. By the definition of $L'$, each $V^{}_j$ is
a subset of one such coset (it is possible that more than one $V^{}_j$
are subsets of a fixed coset). Clearly, we have $V^{}_j\subset y+L'$
for $y\in V^{}_j$. Since $\varphi\in\Phi^{}_{ij}$ is affine, we get
$\varphi(V^{}_j)\subset \varphi(y+L')=\varphi(y)+QL'$. Therefore,
$\Phi[a]$ partitions $\Phi$ into congruence classes induced by the
cosets $L/QL'$. Note that $\Phi[a]\neq\varnothing$ for all $a\in L$.
Note also that even in the one-dimensional case, not all $L^{}_i$ have to be 
equal, cf.~Example \ref{subsec:perseq}. 

The answer to the questions about an LSS being
a model set and being pure point diffractive is now given in the
following theorem. 

\begin{thm}{\cite[Theorem 5.12]{LMS03}} \label{thm:modcoinc} 
  Let $(\V,\Phi)$ (with $\V=(V_1,\ldots, V_m)$) 
  be a primitive LSS, $L'= L_1 + \cdots +
  L_m$, where $L_i:=  \langle V_i - V_i \rangle_{\Z}$. Then the
  following claims are equivalent.  
\begin{enumerate}
\item Each $V_i$ is a regular model set.
\item There is a $k\in \N$ such that a modular coincidence
  relative to $Q^k L'$ occurs in $\Phi_{}^k$.  
\item $\V$ is pure point diffractive (meaning that each
    $V^{}_i$ is pure point diffractive).  \qed
\end{enumerate}
\end{thm} 

Essentially, this theorem tells us that it is sufficient to compute powers
of $\Phi$, until there is a map $\varphi$ which occurs in one row of
$\Phi_{}^k$ only. Alternatively, it follows from \cite{LMS03} that it is
enough to find one sublattice $a+Q^k_{}L'$ which lies entirely in a
single $V^{}_i$, to ensure that each $V^{}_j$ is a regular model set. 
At this point, it is difficult to estimate an upper bound
for $k$, or even to decide whether the derived algorithms will terminate,
assuming that $\V$ is \textit{not} a model set.  
An answer is given by Theorem \ref{thm:upperbd} below.  

Note that the above construction does not yield the
``minimal'' possible internal space, see Example \ref{subsec:abab}. 

\textsc{Remark:} To a multi-component Delone set $\V$ (or a tiling
$\T$), as obtained via a primitive substitution, one can
associate a compact metric space $\mathbb{X}^{}_{\V}$
($\mathbb{X}^{}_{\T}$). Here,
$\mathbb{X}^{}_{\V}:=\cl(\{\V-t\mathbin| t\in\R^d_{}\})$, i.e.,
it is the closure of the set of all translates of $\V$ in the space of 
all multi-component Delone sets, where the metric yields the
following definition of closeness: Two multi-component Delone sets
$\V^{}_1,\V^{}_2$ are close to each 
other if they coincide, after a small shift, inside a large ball around
the origin. This compact space together with the translation action
$\R^d_{}$ yields a \emph{topological dynamical system}
$(\mathbb{X}^{}_{\V},\R^d_{})$. If the eigenfunctions of this action
span a dense subspace of $L^2_{}(\mathbb{X}^{}_{\V})$, then we say
that $\V$ has a \emph{pure point dynamical spectrum}. In the cases
considered here, 
$\V$ (resp.~$\T$) has pure point dynamical dynamical spectrum iff it
is pure point diffractive. This equivalence is often used in the
literature. In fact, this equivalence plays an essential role in the proof of
Theorem \ref{thm:modcoinc}. For details on dynamical systems and this
last connection see~\cite{LMS02,LMS03,Gou05,BL04}.

\section{Main Results}\label{sec:main}

In connection with modular coincidence, we define the following subsets
of $\{1,\ldots ,m\}$:
\begin{equation*}
\Psi^{}_{k}[a]:=\{i\mathbin| (a+Q^k_{}L')\cap V^{}_i \neq \varnothing\}
\quad (a \in L / Q^k L')  
\end{equation*}
The following properties are obvious:
\begin{enumerate}
\item $\Psi^{}_{k}[a]\subset \{1,\ldots,m\}$ and $1\le
  \#(\Psi^{}_{k}[a]) \le m$. 
\item 
  $a+Q^{k}_{}L' \subset \bigcup\limits_{i\in \Psi^{}_{k}[a]}V^{}_i$
\item $(\V,\Phi)$ admits a modular coincidence relative to $Q^k_{}L'$
  iff there is an $a \in L / Q^k L'$ such that $\#(\Psi^{}_{k}[a])=1$. 
\item For $\overline{z}\in L/L'$, denote by $z$ an element of the
  associated coset. Then $\Psi^{}_{0}[z]$ is a partition of
  $\{1,\ldots,m\}$ as $z$ runs through $L/L'$.
\item By $Q^{k+1}_{}L'\subset Q^k_{}L'$ we have
  $\Psi^{}_{k}[a]=\bigcup\limits_{z\equiv a\bmod Q^k_{}L'}
  \Psi^{}_{k+1}[z]$.\\
\end{enumerate}

Already these properties can sometimes be used directly to decide on
the (non-)existence of modular coincidences, see
Example~\ref{subsec:thuemorse}. 

Note that the sets $\Psi^{}_{k}[a]$ can be calculated iteratively:
\begin{align*}
\Psi^{}_{k+1}[b] 
& = \{ i \mathbin| (b+Q_{}^{k+1} L')\cap V^{}_i \neq \varnothing \} \\
& \hspace*{-.9em}\stackrel{\textnormal{Def.~\ref{def:sub_sys}}}{=}
\left\{ i \;\left| 
    \; \left( b+Q_{}^{k+1} L' \right)\cap \left( \bigcup_j
      \bigcup_{\varphi\in\Phi^{}_{ij}[b]} \varphi(V^{}_j) \right)
  \right. \neq \varnothing \right\} \\ 
& = \left\{ i \;\left| \; \bigcup_j
    \bigcup_{\varphi\in\Phi^{}_{ij}[b]} \left( \left( b+Q_{}^{k+1} L'
    \right)\cap \varphi(V^{}_j) \right) \right. \neq \varnothing
\right\} \\  
& \hspace*{-.3em} = \left\{ i
  \;\left| \; \bigcup_j 
    \bigcup_{\varphi\in\Phi^{}_{ij}[b]} \varphi\left( \left( a+Q_{}^k
        L'\right) \cap V^{}_j \right) \right. \neq \varnothing,
  \quad\textnormal{where}\quad\varphi(a)\equiv b \bmod Q_{}^{k+1}L'
\right\} \\  
& = \{ i \mathbin| \textnormal{there exist } j,\;
  \varphi\in\Phi^{}_{ij}[b] \textnormal{ and } a\in L \textnormal{
    s.t. } \varphi(a)\equiv b \bmod Q^{k+1}_{}L' \textnormal{ and
} j\in\Psi^{}_k[a] \},
\end{align*}
where the fourth equation is due to the following: If
$(b+Q_{}^{k+1}L')\cap\varphi(V^{}_j)\neq\varnothing$, then 
there exists an $a$ such that $\varphi(a)\equiv b \bmod
Q_{}^{k+1}L'$. 

For a fixed $\Psi^{}_k[a]$, we call each $\Psi^{}_{k+1}[b]$ which is
determined by the previous equation a \emph{child} of
$\Psi^{}_k[a]$. Conversely, for a fixed $\Psi^{}_{k+1}[b]$, we call all
$\Psi^{}_k[a]$ which appear in the previous equation the
\emph{parents} of $\Psi^{}_{k+1}[b]$.

Example \ref{subsec:kol24} shall illustrate this concept. This 
example also inspires the following notation: We fix a system of
representatives (which includes $0$) in $L$ of $L/QL$, denoted by
$\alpha^{}_i$. Then we can code every $a\in L$ by a finite sequence
$a=(\alpha^{}_1,\alpha^{}_2,\ldots,\alpha^{}_n)$ (the \emph{$Q$-adic
expansion} of $a$), where $a\in \alpha^{}_1+Q L$,
$a-\alpha^{}_1 \in \alpha^{}_2+Q_{}^2 L$, $a-\alpha^{}_1
-Q\alpha^{}_2\in \alpha_3+Q^3_{}L$, etc. 

\begin{defi} \label{def:admiss}
A primitive LSS $(\V,\Phi)$ is called
\emph{admissible}, if all mappings in $\Phi$ have translational parts
of the form $(\alpha^{}_{ijk},0,0,\ldots)$, where $\alpha_{ijk}$ is an
element of a fixed system of representatives for $L/QL$.   
\end{defi}
Note that, by the
definition of an LSS, all such translational
parts appear exactly once in each column of $\Phi$ in the admissible
case. The admissible primitive LSS are those
to be considered in the sequel. Later we shall discuss 
nonadmissible substitution systems (Theorem
\ref{thm:vmldv}). For an example of a nonadmissible primitive LSS, 
cf.~Example \ref{subsec:nonadmiss}.

In the following, we introduce a simple method for computing modular
coincidences in the admissible case.

\begin{defi} \label{def:coincgraph}
The \emph{coincidence graph} $G_{(\V,\Phi)}$ of a primitive
LSS $(\V,\Phi)$ is a directed
\mbox{(multi-)}graph, where the vertices represent sets $\Psi^{}_k[a]$.
If $\Psi^{}_k[a]$ and $\Psi^{}_{\ell}[b]$ are equal as sets, then both are
represented by the same vertex. \\
We mark the $[L:L']$ vertices $\Psi^{}_0[a]$ as ``base points'' (e.g., these
vertices are represented by a double circle instead of an ordinary
circle). We draw an edge from $v_1$ to $v_2$ iff $v_2$ is a child
of $v_1$. Therefore, there are at most
$|\det Q|$ outgoing edges at each vertex. 
\end{defi}

Of course, one can also label the edges by the element $z\in L'/QL'$
which induces the parent/child relation. This is useful if one
wants to analyze the substitution further, e.g., determine the
windows of a model set explicitly. In this case, multiple edges are
possible, and there are exactly $|\det Q|$ outgoing edges at each
vertex. For examples of coincidence graphs, see
Examples \ref{subsec:chairtable}, \ref{subsec:haus}.  

The practical construction of a coincidence graph works as follows. 
Let us introduce the following notation: For an (admissible) LSS
$(\V,\Phi)$ and fixed $j$, denote by
$\Phi(j)^{}_z$ the letter 
$i$ for which $Qx+z\in\Phi^{}_{ij}$ (where $z\in L/QL'$). By definition,
this is (for every $j$ and $z$) exactly one letter. Then we obtain the
coincidence graph as follows: We start with the base points
$\Psi^{}_0[a]$ ($a\in L/L'$). Suppose
$\{s^{}_1,\ldots,s^{}_{\ell}\}=\Psi^{}_0[a]$ for some $a$, then there
is a directed edge to the vertex representing
$\{\Phi(s^{}_1)^{}_z,\ldots\Phi(s^{}_{\ell})^{}_z\}$ for each
$z$. From these new vertices, we proceed further (to the vertices
$\{\Phi_{}^2(s^{}_1)^{}_z,\ldots\Phi_{}^2(s^{}_{\ell})^{}_z\}$) until
no further edges and vertices have to be added. This graph is then the
coincidence graph.

\begin{lemma}\label{lem:children}
Let $(\V,\Phi)$ be an admissible primitive substitution system.
If $\Psi^{}_k[a]=\Psi^{}_{\ell}[c]$, then $\Psi^{}_k[a]$ and
$\Psi^{}_{\ell}[c]$ have the same set of children.
\end{lemma}  

\begin{proof}
Case $L=L'$: Here $a$ in $\Psi^{}_k[a]$ is actually a representative
of $L/Q^{k}_{}L$, wherefore we can write
$(\alpha^{}_{a_1},\ldots,\alpha^{}_{a_k})$ for $a$. Similarly, we can write
$(\alpha^{}_{c_1},\ldots,\alpha^{}_{c_{\ell}})$ for $c$. An (admissible) map
$\varphi(x)=Q x+(\alpha^{}_{ijk},0,0,\dots)$ yields (with these
identifications) 
$\varphi(a)=(\alpha^{}_{ijk},\alpha^{}_{a_1},\ldots,\alpha^{}_{a_k})$
and
$\varphi(c)=(\alpha^{}_{ijk},\alpha^{}_{c_1},\ldots,\alpha^{}_{c_{\ell}})$. 
Now, since $\Psi^{}_k[a]=\Psi^{}_{\ell}[c]$, the same maps operate on
$a+Q^k_{}L$ as well as on $c+Q^{\ell}_{}L$, and we get
$\Psi^{}_{k+1}[t+Qa]=\Psi^{}_{\ell+1}[t+Qc]$ (where $t$ is an associated
translational part). Furthermore, we observe in the admissible case
that each child has exactly one parent
($\Psi^{}_{k+1}[(\alpha^{}_1,\ldots,\alpha^{}_{k+1})]$ has parent
$\Psi^{}_{k}[(\alpha^{}_2,\ldots,\alpha^{}_{k+1})]$). This proves the
lemma for $L=L'$.

Case $L\neq L'$: Since each set $\Psi^k_{}[a]$ determines the
representative of $L/L'$ uniquely, we can simply index the position by
the $Q$-adic expansion relative to $L'$. This reduces this case to the
previous one.
\end{proof} 

It follows that we need to compute children only as long as we find
sets that did not occur before. This justifies the representation of
equal children $\Psi^{}_k[a]=\Psi^{}_{\ell}[b]$ by a single vertex
in Def.~\ref{def:coincgraph}. As a consequence, we obtain the
following result. 

\begin{thm} \label{thm:aequiv}
Let $(\V,\Phi)$ be an admissible primitive LSS. Then,
$(\V,\Phi)$ admits a modular coincidence iff at least one
vertex of the coincidence graph $G_{(\V,\Phi)}$
represents a singleton set. The length of the shortest path from
one of the vertices $\Psi^{}_0[a]$ to this singleton set is the minimal
$k$ such that a modular coincidence relative to $Q^kL'$ occurs
in $\Phi_{}^k$. \hfill $\square$ 
\end{thm}

The claim follows immediately from Lemma \ref{lem:children} and the
construction of $\Psi$. This result yields an upper bound for the
number of iterations needed to compute modular coincidences. 

\begin{thm} \label{thm:upperbd}
Let $(\V,\Phi)$ be an admissible primitive LSS with $m \ge 2$ components.
If $\Phi^k_{}$ admits no modular coincidence for $k \le 2^m-m-2$,
then $\Phi^k_{}$ admits no modular coincidence for any $k \in \N$. 
\end{thm}

\begin{proof} 
As described above, it is possible to determine the sets $\Psi^{}_k[a]$
iteratively. The LSS $(\bs{V},\Phi^k_{})$ admits a modular coincidence
iff there is an $a$ such that $\Psi^{}_k[a]$ is a singleton set. 

From Lemma~\ref{lem:children}, it follows that we need to compute
children until no new sets appear. 
Since every child is a nonempty subset of 
$\{1,2, \ldots, m\}$, after at most $2^m-2$ iterations all possible
sets have occurred. (We count the determination of $\Psi^{}_0[\cdot]$ as
the zeroth iteration step. So $\Psi^{}_1[\cdot]$ corresponds to
$\Phi_{}^1$.) Moreover, after at most $2^m-m-2$, iterations either a
set with just one element occurred, or, by Lemma \ref{lem:children},
it will never occur. By Theorem \ref{thm:aequiv}, this is
equivalent to the claim of the theorem. 
\end{proof}  

As a consequence of Theorem \ref{thm:aequiv}, we obtain some
conditions for $(\V,\Phi)$ admitting a modular coincidence. By Theorem
\ref{thm:modcoinc}, this translates to necessary conditions for $\V$
consisting of model sets, resp.~being pure point diffractive.    

\begin{corollary} \label{cor:psising} 
Let $(\V,\Phi)$ be an admissible primitive LSS. 
If there is an $a$ such that $\Psi^{}_0[a]$ is a singleton set,
then $\V$ consists of model sets. 
\end{corollary} 
\begin{proof}
This is an immediate consequence of Theorem \ref{thm:aequiv}: 
$\Psi^{}_0[a]$ is already a singleton set, thus, by the
construction of $\Psi$, we obtain $V_i=a+L'$ for some $1 \le i \le
m$. By the remark  following Theorem \ref{thm:modcoinc}, $\V$ consists
of model sets.   
\end{proof}

\begin{corollary} \label{cor:paircoinc}
Let $(\V,\Phi)$ be an admissible primitive LSS,
such that no $\Psi^{}_0[a]$ is a singleton set.
If $\V$ consists of model sets, then $\Phi$ admits at least one
\emph{pairwise} coincidence in some $\Psi^{}_0 [a]$, i.e., there are 
$j \ne k$ such that $j,k \in \Psi^{}_0[a]$ for some $a$ and 
$\Phi_{ij} \cap \Phi_{ik} \ne \varnothing$ for some $i$.
\end{corollary}
\begin{proof} 
Assume that $\Phi$ admits no pairwise coincidence in any $\Psi^{}_0[a]$.
Therefore, and because any set $\Psi^{}_0[a]$ contains more than one element,
any child $\Psi^{}_1[Qa+z]$ contains more than one
element. Inductively, no set $\Psi^{}_k[\cdot]$ is a singleton set. 
\end{proof}

Surprisingly, the last corollary applies to a large number of
substitutions found in the literature, showing that they do \textit{not}
consist of model sets, since they admit no pairwise coincidence
(cf.~Examples \ref{subsec:thuemorse},
\ref{subsec:chairtable}, \ref{subsec:natalies}).  
In \cite[Def.~1.2]{PFpre}, the author defines \emph{bijective substitutions}.
This definition reads in our notation as follows. 

\begin{defi}
Let $(\V,\Phi)$ be an admissible primitive LSS.
If each $\varphi$ occurring in $\Phi$ appears exactly once in
each row, then $(\V,\Phi)$ is called \emph{bijective substitution
  (system)}.  
\end{defi}

The following statement is mentioned in \cite[Section 3.3]{PFpre}, but
without proof. Here we obtain the result quite easily.

\begin{corollary} \label{cor:bij}
Let $(\V,\Phi)$ be an admissible primitive LSS,
such that $\Phi$ is bijective and $\V$ is nonperiodic. 
Then $\V$ does not
consist of model sets, therefore $\V$ is not pure point diffractive. 
\end{corollary}
\begin{proof}
The proof is an application of Corollary \ref{cor:paircoinc}. 
$\Phi$ is bijective, thus $\Psi$ admits no pairwise coincidence at
all. All we need to show is that no $\Psi^{}_0[a]$ is a singleton set.  

First we note that in $\V$ all types of points occur with the same
frequency, if $\V$ arose from a bijective substitution:  
Since $(\V,\Phi)$ is bijective, every $\varphi$ appears exactly once
in each row of $\Phi$. Since $(\V,\Phi)$ is admissible, every
$\varphi$ appears exactly once in each column of $\Phi$. 
Therefore, the row sums of the substitution matrix 
$S(\Phi)=(\# \Phi_{ij})_{1 \le i,j   \le m}$  are all equal. Let 
$q \in \N$ be this row sum. Thus $q$ is an eigenvalue of $S(\Phi)$,
with corresponding normalized eigenvector $(1/m,1/m,\ldots,1/m)^T$.
Since $S(\Phi)$ is primitive, it follows from the Perron-Frobenius
Theorem that $q$ is the unique eigenvalue of $S(\Phi)$ which is
largest in modulus (since it is the unique one with a strictly positive
eigenvector). It is well-known (and easy to see) that the
corresponding (normalized) eigenvector contains the relative
frequencies of the points of each type. Thus, all types of points have
the same frequency in $\V$. 

Assume that there is some $\Psi^{}_0[a]$ such that $\# \Psi^{}_0[a]=1$. This
means that all points and only the points in $a+L'$ are of the same
type, which occurs with relative frequency $\frac1m$. Therefore, the
index $[L:L']$ equals $m$, but there are exactly $m$ types, and each coset
$L/L'$ contains at least one type. Then the only possibility is $\#
\Psi^{}_0[b]=1$ for all $b$ and $ \V$ is periodic (w.r.t. the lattice $L'$). 

Under the nonperiodicity assumption of the corollary, it now follows that each
$\Psi^{}_0[a]$ contains more than one element, which completes the proof.  
\end{proof}

The following theorem allows us to extend the results to nonadmissible 
LSS. Under a certain condition, we can find for a
primitive LSS $(\V,\Phi)$ an admissible one,
say, $(\V',\Phi')$, such that $\V$ and $\V'$ are MLD. 
Here and in the sequel, $\V \cap R\B^d$ is to be read as 
$(V_1 \cap R\B^d, \ldots, V_m \cap R\B^d)$.

\begin{defi} \label{def:nicegrow}
Let $(\V,\Phi)$ be a primitive LSS 
and $F$ a fundamental domain of the lattice $\supp(\V)$. We say that 
$(\V,\Phi)$ is \emph{nicely growing}, if there is an $R>0$ such that, 
for all $R$--clusters $\V \cap (t + R \B^d)$,
\begin{equation} \label{eqn:nicegrow} 
\supp(\V) \cap (Qt + QF + R \B^d)  \subset \supp 
\left( \Phi ( \V \cap (t + R \B^d )) \right). 
\end{equation}  
\end{defi}  
In plain words: In particular, we require that the action of $\Phi$
maps the set of points of $\V \cap R \B^d$ to the points inside a
region which is considerably larger, namely, $\V \cap (QF + R
\B^d)$ (and maybe to some more points).  
And in general, the points of $\V$ inside any ball of radius $R$, centered
at $t$, are mapped to the points in the set $\V \cap (Qt + QF + R
\B^d)$, ``centered'' at $Qt$ (and maybe to some more points).  

\textsc{Remark:} All admissible substitution systems are nicely
growing. If $(\V,\Phi)$ is admissible, then $\Phi$ maps $\V \cap R\B^d$ 
onto $\V \cap Q(R \B^d)$ (and possibly more), at least for $R$ large
enough. (E.g., the admissible substitution of the substitution system
$(\V,\Phi)$ in Example \ref{subsec:perseq} maps the cluster 
$\V \cap  \B^1 = \V \cap [-1,1]$ to the cluster 
$\V \cap [-6,11]$, which is a superset of 
$\V \cap ([0,6] + [-1,1]) = \V \cap [-1,7]$,
thus \eqref{eqn:nicegrow} is fulfilled in this example for $t=0$,
$R=1$.) Then the number of points in the support of the left hand 
side of \eqref{eqn:nicegrow} grows approximately like $R^d+cR^{d-1}$, 
where the number of points in the support of the 
right hand side grows like $R^d |\det(Q)|$.  
Since $|\det(Q)| >1$, Equation \eqref{eqn:nicegrow} is fulfilled for
$R$ large enough. This shows that all admissible
substitutions are nicely growing.
  
This is also true for the nonadmissible substitutions in Example
\ref{subsec:nonadmiss}. Essentially, \textit{all} primitive lattice
substitution systems known to the authors are nicely growing (but
cf.~Example \ref{subsec:notnice}). Note that the requirement of
checking $R$--clusters is a finite problem, since there are
only finitely many translation classes of them.

\begin{thm}\label{thm:MLDtiling} 
Let $(\V,\Phi)$ be a nicely growing primitive LSS.
Then  there exists $(\V',\Phi')$, such that $\V$ and $\V'$ are
MLD and $(\V',\Phi')$ is admissible. 
\end{thm} 

Note that our general assumption on $\V$ --- all 
$r$--clusters in $\V$ are contained in some iteration  
$\Phi^k({\bs v})$ of a single point  ${\bs v}$ ---  ensures that $\V$
is repetitive. 

\begin{proof}
The basic idea is to construct $\V'$ out of $\V$ by assigning to each 
colored point in $\V$ a colored point in $\V'$ at the same position,
where the color of the latter one depends on the local surrounding of
the former one. This construction is well--known, see for example~\cite{kp}
(``forcing the border'') or \cite{ors} (``collaring''). In general, $\V'$ will
use more colors than $\V$. Then the admissible MFS $\Phi'$ arises in
a quite natural way, using the property $\Phi(\V)=\V$. 

Let $R>0$ such that 
\begin{equation*}
\supp(\V) \cap (QF + R \B^d) \subset
\supp  \left( \Phi ( \V \cap R \B^d ) \right).
\end{equation*}
Since $(\V,\Phi)$ is nicely growing, we find such an $R$ by
Definition \ref{def:nicegrow}. Let $u_1,\ldots,u_d$ be a basis of the
lattice $L$, i.e.,  $L:= \langle u_1, \ldots, u_d \rangle_{\Z}$. We
introduce the fundamental domain 
\begin{equation*}
F=\{ \sum_{i=1}^d \lambda_i u_i \mathbin| 0 \le \lambda_i < 1 \}
\end{equation*}
of $L$. For a fixed lattice $L$, the lattice basis is not
unique, hence there are many possibilities how $F$ may look like. 
For simplicity, we may choose the vectors $u_i$ such that the diameter
$\diam(F)$ becomes minimal among all possible fundamental domains. 
But this is not relevant in the following. 

Let $x\in L$. Then any set $x+F$ contains exactly one element of $L$,
namely $x$, and any set $Qx+QF$ contains exactly $|\det(Q)|$ elements of
$L$, which just are representatives of $L/QL$. 

Consider all types of $R$--clusters in $\V$, i.e., all sets 
\begin{equation} \label{eq:vcapb}
\V \cap (x+ R \B^d) := (V_1 \cap (x+ R \B^d), \ldots , V_m \cap 
(x+ R  \B^d)), 
\end{equation}
where $x \in L$. There are finitely many
translation classes of these, $[C^{(1)}]$, $[C^{(2)}]$, $\ldots$,
$[C^{(n)}]$, say.  For simplicity, we represent every translation class by a
cluster $C^{(i)}$ centered at the origin. Note that the $C^{(i)}$ are of the
form $(C^{(i)}_1, \ldots, C^{(i)}_m)$. Now we construct $\V'$ out of
$\V$ in the following way:   
\begin{equation*}
V'_i:= \{ x \in L \mathbin| \V \cap (x+ R \B^d) \in [C^{(i)}] \}, 
\end{equation*} 
and $\V'=(V'_1, \ldots, V'_n)$.  
In plain words: If $x$ is the center of a cluster of type $C^{(i)}$ in
$\V$, then $x$ in $\V'$ gets color $i$. Since all lattice points $x
\in L$ are covered, we get $L = \bigcup_{1 \le i \le n}
V'_i$. Clearly, any $x \in L$ is the center of exactly one
$R$--cluster, therefore this union is disjoint. 
 
The MFS $\Phi'$ is given as follows: Let $\V' \cap \{ x \}$ be of type
$j$, i.e., $x \in V'_j$. Now define $\Phi'$ by 
\begin{equation*}
\Phi' (\V' \cap \{x\}) := \V \cap (Qx+QF) = ((Qx +QF) \cap V_1',
\ldots,(Qx +QF) \cap V_n').
\end{equation*} 
The sets $(Qx +QF) \cap V_i'$ contain finitely many points, which 
define the maps in $\Phi'_{ij}$ as follows. Let $j$ be fixed. 
If $(Qx +QF) \cap V_i' = \{ b_1 , b_2 , \ldots, b_k\}$, 
then set $a_{ij\ell} := b_{\ell} - Qx \in QF\cap L$,
$\varphi_{\ell}(y)=Qy + a_{ij\ell}$ ($1 \le \ell \le k$) and 
$\Phi'_{ij} = \{\varphi_1,\varphi_2, \ldots, \varphi_k\}$.

We claim that $(\V',\Phi')$ is an admissible primitive LSS
such that $\V$ and $\V'$ are MLD. The maps of $\Phi'$ are of the form  
$x \mapsto Qx+a_{ij\ell}$. By construction, 
the substitution given by $\Phi'$ maps every element of $\V'$ to
$|\det(Q)|$ elements, which are contained in some translate of $QF$.
Moreover, every $a_{ij\ell}$ is an element of $QF\cap L$, which is a
fixed system of representatives of $L/QL$. Thus the definition of
admissibility is fulfilled. 

Obviously, a point in $\V'$ determines an $R$--cluster in $\V$, so
$\V$ is locally derivable from $\V'$. By construction, $\V'$ is
locally derivable from $\V$, so $\V$ and $\V'$ are MLD.

Let $\V \cap \{x\}$ be of type $i$, $R>0$. Since $(\V,\Phi)$ is
primitive, and by our general assumption (see the beginning of Section
\ref{sec:modsets}), $\Phi_{}^k(\V \cap \{x\})$ contains
all types of $R$--clusters for some appropriate $k$. 
Therefore, applying $(\Phi')^k$ to the
corresponding point $\V' \cap \{x\}$ yields all types of points in
$\V'$. The choice of $\V \cap \{x\}$ is arbitrary, so it can
be the center of any type of $R$--cluster in $\V$. Altogether, this
means that, for any $\V' \cap \{x\}$, there is $k \in \N$ such that
$(\Phi')^k (\V' \cap \{x\})$ contains all types of points that occur in
$\V'$. This is equivalent to $\Phi'$ being primitive. 

All we are left with to show is that $\Phi'$ is well-defined, i.e.,
that two points of the same type are substituted in the same
way. This follows essentially from the condition $\Phi(\V)=\V$: Since
two congruent clusters in $\V$ are mapped by $\Phi$ to congruent
clusters  in $\Phi(\V)$, the corresponding points
in $\V'$ are mapped to congruent clusters in $\V'$. Formally, we need
to show: If $x,y \in V'_i$ for some $i$, i.e., 
$\V' \cap \{x\} =  (\V' \cap \{y\})  + x-y$, then 
\begin{equation} \label{eq:mfseind}  
 \Phi' (\V' \cap \{x\}) =  \Phi' (\V' \cap \{y\})  + Q(x-y) 
\end{equation}
holds. So, let $\V' \cap \{x\} =  (\V' \cap \{y\})  + x-y$. By
construction of $\V'$, it follows that the corresponding $R$--clusters
in $\V$ are of the same type, i.e.,
\begin{equation*}
\V \cap (x+R\B^d) = (\V \cap (y+R\B^d)) + x-y.
\end{equation*} 
Since $\Phi(\V)=\V$, we obtain
\begin{equation*}
\Phi \left( \V \cap (x+R\B^d) \right) = \Phi \left( (\V \cap
  (y+R\B^d)) + x-y \right) = \Phi \left( \V \cap
  (y+R\B^d) \right) + Q(x-y).
\end{equation*}
Consider $\V \cap (Qx+QF+R\B^d), \V \cap (Qy+QF+R\B^d) + Q(x-y)$. The
growing nicely condition implies:
\begin{equation*}
\begin{array}{lcl} 
\supp \left( \V \cap (Qx + QF+R\B^d) \right) & \subset &  \supp
\left( \Phi \left( \V \cap (x+R\B^d) \right) \right) \\ 
 & & \qquad \| \\  
\supp \left( \V \cap (Qy + QF+R\B^d) \right) +Q(x-y)  &\subset & \supp
\left( \Phi \left( \V \cap   (y+R\B^d) \right)\right)  +Q(x-y).  \\
\end{array} 
\end{equation*}
The four sets $\V \cap (Qx + QF+R\B^d)$, $\Phi \left( \V \cap
  (x+R\B^d) \right)$, $\V \cap (Qy + QF+R\B^d)$ and 
$\Phi \left( \V \cap  (y+R\B^d) \right)$  
occurring in the above equation are all subsets
of $\V$ as multi-component sets, i.e., with respect to colors.   
We obtain
\begin{equation*}
\begin{array}{lcl} 
\V \cap (Qx + QF+R\B^d) & \subset &  \Phi \left( \V \cap (x+R\B^d)
\right) \\  
 & & \qquad \| \\  
\V \cap (Qy + QF+R\B^d) +Q(x-y)  &\subset & 
\Phi \left( \V \cap   (y+R\B^d) \right)  +Q(x-y).  \\
\end{array}
\end{equation*}
Obviously, the supports of the left hand sides are equal. Therefore,
it follows from the above equation, that the left hand sides are equal as
multi-component sets: 
\begin{equation*}
\V \cap (Qx + QF+R\B^d)  = \V \cap (Qy + QF+R\B^d) + Q(x-y)
\end{equation*}
By the construction of $\V'$, the $R$--cluster $C_R = \V \cap
(x+R\B^d)$ in $\V$ determines the color of $x \in \V'$. Therefore,
the above equation yields 
\begin{equation*}
\V' \cap (Qx + QF) = \V' \cap (Qy + QF) +Q(x-y),
\end{equation*} 
and by the construction of $\Phi'$, this is equivalent to
(\ref{eq:mfseind}). Altogether, $(\V',\Phi')$ is an admissible
primitive LSS, such that $\V$ and $\V'$ are
MLD.  
\end{proof}

The last result proves useful as follows. 

\begin{thm} \label{thm:vmldv}
Let $(\V,\Phi)$ be a primitive LSS, and
$(\V',\Phi')$ be an admissible one, derived from $(\V,\Phi)$ as in the
proof of Theorem {\rm \ref{thm:MLDtiling}}. If $\V'$ consists of model sets,
then so does $\V$.   
\end{thm}
\begin{proof}
Let $\V=(V_1,\ldots,V_m)$. By the construction of $\V'$, each set $V_i$ is
partitioned into sets $V'_{i_1}, V'_{i_2}, \ldots, V'_{i_k}$ belonging
to $\V'$. Recall: An element of $\V'$ at position $x$ got its type
according to the $R$--cluster surrounding $x$ in $\V$. 
Denote all translation classes of $R$--clusters with center of
type $i$ by $[C_{i_1}], [C_{i_2}], \ldots, [C_{i_k}]$. Then $V_i =
V'_{i_1} \cup V'_{i_2} \cup \cdots \cup V'_{i_k}$. 

Let $\V'$ consist of model sets. This means that each $V'_{i_j}$ is a
model set. To each $V'_{i_j}$ belongs a window set $W'_{i_j}$. Let $V'_{i_1},
V'_{i_2}, \ldots, V'_{i_k}$ be the sets assigned to $V_i$, i.e., 
$V_i = \bigcup_{1\le j \le k} V'_{i_j}$. Let $\Omega'_{i_j}$ be the
window of $V'_{i_j}$ (cf.~\ref{def:modset}). Then
\begin{equation*}
V_i =  \bigcup_{1\le j \le k} V'_{i_j} = \bigcup_{1\le j \le k}
\varLambda(\Omega'_{i_j}) = \varLambda(\bigcup_{1\le j \le k}
\Omega'_{i_j}).
\end{equation*} 
Since the windows $\Omega'_{i_j}$ have nonempty interior and compact
closure, this is also true for their (finite) union. Thus $V_i$ is a
model set, too. 
\end{proof}

\section{Remarks on Complexity}\label{sec:complex}

Theorem \ref{thm:upperbd} gives an exponential upper bound for computing
modular coincidences. This may be far away from the best upper bound. 
In particular, it is not known whether the required number of iterations is
polynomial in $m$ or not. 

If $k=2^m$, the $k$--th power of an $(m \times m)$--MFS $\Phi$ can be
computed by $m-1$ ``matrix multiplications'', to be understood as stated
after Def.~\ref{def:mfs}. If $k' < k$ and $\Phi_{}^{k'}$ already admits a
modular coincidence, then $\Phi_{}^k$ also does. Therefore, it suffices to
compute $\Phi, \Phi_{}^2, \Phi_{}^4, \Phi_{}^8 \ldots$. At first sight, this
would lead to a polynomial bound for computing modular coincidences.   
However, unlike usual matrix multiplications, the costs for each
multiplication is not $m^3$, but the costs per multiplication increase
with $k$. We do not study this question in detail here, but
essentially, the costs grow faster than the sum of the entries of the
substitution matrix $S(\Phi_{}^k)$.   

The authors conjecture that the lowest upper bound is quadratic in
$m$, hence not linear in $m$. Consider the following example:
\begin{equation} \label{eq:worstex} 
  1 \mapsto 1\; 2, \quad 2 \mapsto 2\; 3, \quad 3 \mapsto 3\;4,
 \; \ldots, \; m-1 \mapsto m-1 \; m , \quad m \mapsto 1\;1 
\end{equation}
If $m=2$ this reads as the well-known substitution $a \to ab, \; b \to
aa$. For $2 \le m \le 16$, the first modular coincidence occurs in
$\Phi_{}^{(m-1)^2}$. This is a hint that a general upper bound is
quadratic in $m$. It will be difficult to find a worse example, since
the dimension (=1) and the factor $\det(Q)$ (=2) are the smallest
possible for admissible primitive LSS, as
well as the number of pairwise coincidences (=1). For 
$m \in \{2,3,4,5,6\}$, we computed for {\em all} substitutions with these
three properties the minimal number $k$, such that $\Phi_{}^k$ shows a modular
coincidence (if there is such a $k$). It turned out that the
substitution (\ref{eq:worstex}) yields the maximum value of $k$ in
these cases.

\section{Examples of Sequences and Tilings}\label{sec:examples}

\subsection{Periodic substitution I}\label{subsec:perseq}
The primitive substitution given by the following MFS yields a periodic
sequence. Here $h_i$ denotes the map $h_i(x)=6x+i$, $i \in \Z$. 
\begin{equation*}
\Phi = 
\begin{pmatrix} \{h_0,h_3\} &  \{h_0,h_3\} &  \{h_0,h_3\} \\
\{h_1, h_5\} & \{h_1, h_5\} & \{h_1, h_5\} \\
\{h_2,h_4\} & \{h_2,h_4\} & \{h_2,h_4\} \\
\end{pmatrix}
\end{equation*}
In symbolic notation, this reads $a\mapsto abcacb$, $b\mapsto abcacb$,
$c\mapsto abcacb$. This substitution on $L=\Z$ yields 
$L^{}_a=3\Z$, $L^{}_b=2\Z$ and $L^{}_c=2\Z$ (and therefore
$L'=L$). Note that here --- as in some of the following examples ---
we use letters $a$, $b$, $c$ instead of numbers $1$, $2$, $3$.
To obtain a tiling, replace the letter at position $n$ by the interval
$[n,n+1]$ (and label it with its letter). To obtain a multi-component
Delone set, let $\V=(V_a,V_b,V_c)$, where $V_a$ contains all positions
of letters $a$ etc. 

\subsection{Periodic substitution II}\label{subsec:abab}
The trivial substitution $a\mapsto ab$, $b\mapsto ab$ on $L=\Z$ yields the
periodic solution $\ldots ababababababa\ldots$ (and therefore
$V^{}_a=2\Z$, $V^{}_b=2\Z+1$ and $L'=2\Z$). It already shows a modular
coincidence relative to $QL'$ (each $\Phi[i]$, $0\le i\le 3$, consists
of exactly one element). Following the construction preceding
Def.~\ref{def:modcoinc}, one would obtain as internal space
$H=C^{}_2\times\Z^{}_2$ (where $C^{}_2=\Z/2\Z$ is the cyclic group of
order $2$, and $\Z^{}_2$ are the $2$-adic integers) and the lattice
$\tilde{L}=\{(t,t\bmod 2,t)\in \R\times C^{}_2\times\Z^{}_2\mathbin|
t\in\Z\}$, but it is easy to see that $H'=C^{}_2$ and
$\tilde{L}'=\{(t,t\bmod 2)\in\R\times C^{}_2\mathbin| t \in \Z \}$
suffices.

\subsection{Thue-Morse sequence}\label{subsec:thuemorse}
The Thue-Morse sequence is defined as the fixed point
beginning by $0$ of the substitution $0\mapsto 01$, $1\mapsto 10$
(it starts $0110100110010110100\ldots$). There is also an arithmetic
definition of this sequence (e.g., see~\cite[p.\ 2]{Haeseler},
\cite[2.1.1]{Fogg}): For any integer $n$, denote by $S^{}_2(n)$ the
sum of the dyadic (binary) digits of $n$, i.e., $S^{}_2(n)=\sum_{i\ge
  0}n^{}_i$, if $n=\sum_{i\ge 0}n^{}_i 2^i_{}$
($n^{}_i\in\{0,1\}$). Then the $n$-th symbol in the sequence is a $1$
iff $S^{}_2(n)$ is odd. Obviously, for $n<2^k_{}$, we have the
following: if $S^{}_2(n)$ is odd, then $S^{}_2(n+2^k_{})$ is even and
vice versa. Therefore $\Psi^{}_k[n]=\{i\mathbin| (n+2^k_{}\Z)\cap
V^{}_i\}=\{0,1\}$ for all $k$ and $n$, and the Thue-Morse substitution
admits no modular coincidence. 

\subsection{Parents and children}\label{subsec:kol24}
Consider a substitution defined by $a\mapsto aba$,
$b\mapsto bcc$ and $c\mapsto abc$ (this substitution arises in the
study of the Kolakoski sequence over the alphabet $\{2,4\}$,
see~\cite{Sin03}). The corresponding MFS $\Phi$ is given by:
\begin{equation*}
\Phi= \begin{pmatrix}
\{h_0,h_2\} & \varnothing & \{h_0\} \\
\{h_1\} & \{h_0\} & \{h_1\} \\
\varnothing & \{h_1,h_2\} & \{h_2\} \\
\end{pmatrix},
\end{equation*}
with $h_i(x)=3x+i$. A fixed 
point is given by $\ldots aba\, bcc\, abc\, \dot{a}ba\, bcc\, aba
\ldots$, where the dot marks the zeroth position. It is easy to check
that here we have $L=\Z=L'$. In this case, as $L=L'$ and $\Phi$ is
admissible, the check for modular coincidence reduces to finding a map
$g$ which is shows up in only one row of $\Phi^k$ for some $k \in
\N$. Therefore we get (note that
$\Phi^{}_{aa}=\{h_0,h_2\}$, $\Phi^{}_{ba}=\{h_1\}$, etc.)
\begin{equation*}
\Phi[0]\subset \Phi^{}_{aa}\cup\Phi^{}_{bb}\cup\Phi^{}_{ac},\quad
\Phi[1]\subset \Phi^{}_{ba}\cup\Phi^{}_{cb}\cup\Phi^{}_{bc},\quad
\Phi[2]\subset \Phi^{}_{aa}\cup\Phi^{}_{cb}\cup\Phi^{}_{cc}.
\end{equation*}
Here we use letters rather than numbers to avoid confusion, because 
numbers are used to represent cosets below. Since none of the sets
$\Phi[\cdot]$ is contained in just one row of the MFS, $\Phi$ itself
does not show a modular coincidence. (For instance, $\Phi[0]$ is
contained in the first and second row of $\Phi$, $\Phi[1]$ in the
second and last row, etc.) Nevertheless, $\Phi^2$ shows a
modular coincidence already, or more precisely, three of them.   
\[ \Phi^2 = \begin{pmatrix} 
\{ h^2_0, h^{}_0  h^{}_2, h^{}_2  h^{}_0 , h^2_2 \} & 
  \{ h^{}_0 h^{}_1 , h^{}_0  h^{}_2 \} & 
  \{ h^2_0, h^{}_2  h^{}_0, h^{}_0  h^{}_2 \} \\ 
\{ h^{}_1 h^{}_0, h^{}_1  h^{}_2, h^{}_0  h^{}_1  \} &
  \{ h^2_0,  h^2_1 , h^{}_1  h^{}_2 \} &
  \{ h^{}_1 h^{}_0, h^{}_0  h^{}_1, h^{}_1  h^{}_2  \} \\
\{ h^2_1, h^{}_2  h^{}_1 \} & 
  \{ h^{}_1 h^{}_0, h^{}_2  h^{}_0, h^{}_2  h^{}_1 , h^2_1 \} & 
  \{ h^2_1, h^{}_2  h^{}_1 , h^2_2 \} 
\end{pmatrix} \]
The map $h^{}_0 \circ h^{}_2 (x) = 9x+6$ is contained in just the first
row, the map $h^{}_1 \circ h^{}_2 (x) = 9x+7$ just in the second row
and the map $h^{}_2 \circ h^{}_1 (x) = 9x+5$ just in the third row. 

Since $L=L'$, we have $\Psi^{}_0[0]=\{a,b,c\}$.
The next step is easy: For $\Psi^{}_{1}[i]$ we just have to consider
the maps $\Phi[i]$ (see the first indices of the matrix elements of
$\Phi$ involved in $\Phi[i]$) to obtain
\begin{equation*}
\Psi^{}_1[0]=\{a,b\},\quad
\Psi^{}_1[1]=\{b,c\},\quad
\Psi^{}_1[2]=\{a,c\}.
\end{equation*}
Basically, we make use of the (unique) ternary expansion of each
natural number and that the maps of the MFS line up with this ternary
expansion. E.g., to calculate $\Psi^{}_2[7]$, we observe that $7=2\cdot
3+1$, therefore $\Psi^{}_2[7]$ is determined by the $i$'s of
$\Phi^{}_{ij}[1]$ with $j\in \Psi^{}_1[2]$, yielding
$\Psi^{}_2[7]=\{b\}$. For completeness, we list all $\Psi^{}_2$:
\begin{equation*}
\begin{split}
& \Psi^{}_2[0]=\{a,b\},\quad
\Psi^{}_2[1]=\{b,c\},\quad
\Psi^{}_2[2]=\{a,c\},\quad
\Psi^{}_2[3]=\{a,b\},\\
& \Psi^{}_2[4]=\{b,c\},\quad
\Psi^{}_2[5]=\{c\},\quad
\Psi^{}_2[6]=\{a\},\quad
\Psi^{}_2[7]=\{b\},\quad
\Psi^{}_2[8]=\{a,c\}.
\end{split}
\end{equation*}
Here, the children of $\Psi^{}_1[i]$ are $\Psi^{}_2[3\cdot i]$,
$\Psi^{}_2[3\cdot i+1]$ and $\Psi^{}_2[3\cdot i+2]$. Since
$\Psi_2[5]$, $\Psi_2[6]$ and $\Psi_2[7]$ are singleton sets, we have a
modular coincidence.  

Note the correspondence between $\Phi^{k}[i]$ and $\Psi^{}_k[j]$. For
instance, the fact that 
$\Psi^{}_2[5]=\{c\}$ corresponds to the fact that $h^{}_2 \circ
h^{}_1 (x) = 9x+5$ is contained in the last (the ``$c$--th'') row
of~$\Phi^2$.

\subsection{Nonadmissible substitution I}\label{subsec:nonadmiss}
Consider the one-dimensional substitution given by the following MFS
(where $h_i$ denotes the map $h_i(x)=3x+i$):
\begin{equation*}
\Phi = 
\begin{pmatrix} \{h_0\} &  \{h_1\} & \{h_2\} \\
\{h_1, h_5\} & \{h_0\} & \{h_1,h_5\} \\
\{h_2\} & \varnothing &  \{h_0\} \\
\end{pmatrix}
\end{equation*}
In symbolic notation, this corresponds to the substitution 
\begin{equation*}
a \to abc \_ \, \_ b, \; b \to ba, \; c \to cba \_ \, \_ b,
\end{equation*}
where \_ stands for an empty position. 
A fixed point of this substitution is  
\begin{equation} \label{eq:nonadmseq}
\ldots cba\, bab\, abc\, bab\, \dot{a}bc\, bab\, cba \, bab \, abc
\ldots 
\end{equation}  
The dot marks the zeroth position. Obviously, this substitution is
not admissible. Nevertheless, a brief look at the first
iterations $\Phi_{}^k(\{0\},\varnothing, \varnothing)$ shows that this
substitution is nicely growing. It is easy to see that this
substitution is equivalent to the admissible substitution
\begin{equation*}
a \to abc, \; b \to bab, \; c \to cba.
\end{equation*}
Both substitutions yield the same sequences, therefore the sequences
are clearly MLD. Here we did not use the construction of
Theorem \ref{thm:vmldv}. This would yield another sequence $\V'$ MLD to
(\ref{eq:nonadmseq}), but it would use more than three letters. Namely,
the construction works with $R=2$. There are eight different types of
$2$--clusters (five-letter words) in (\ref{eq:nonadmseq}), yielding
eight different letters in the new sequence $\V'$. This example is
easily generalized to more than one dimension.

\subsection{Nonadmissible substitution II}\label{subsec:notnice}
The following MFS yields another nonadmissible substitution:
\begin{equation*}
\Phi = 
\begin{pmatrix} 
\{h_{15}\}    & \{h_{15}\} & \varnothing   & \{h_{15}\}    \\
\{h_{14}\}    & \{h_{13}\} & \{h_{12}\}     & \{h_{14},h_{13}\} \\
\{h_{13}\}    & \{h_{14}\} & \{h_{15},h_{14},h_{13}\} & \varnothing   \\
\{h_{12},h_0\}& \{h_{12}\} & \varnothing    & \varnothing
\end{pmatrix}, 
\end{equation*}
where $h_i(x)=4x+i$. In symbolic notation, this would look quite weird,
since the action of $\Phi$ contains large shifts to the right. E.g., a
letter $a$ at position 0 becomes a word $d \, \_ \, \ldots \,
\_ \, dcba$ under the action of $\Phi$, showing a large gap between
the first letter $d$ at position 0 and the second letter $d$ at
position 12. This substitution
is clearly not admissible. Nevertheless, there is a $\V$ such that
$\Phi(\V)=\V$. E.g., let $C=(\{-5\},\varnothing,\{-4\},\varnothing)$.
The iteration of $\Phi$ on $C$ yields a sequence with $\Phi(\V)=\V$,
namely $\V=\bigcup_{k \in \N} \Phi_{}^k(C)$. Thus $(\V,\Phi)$ is a primitive
LSS. It is not easy to see whether this substitution is nicely
growing. If so, it can be shown that $R$ must be at least 15. (Note
that the difference between the smallest and the largest
translational part occurring in $\Phi$ is 15.)

This example can be generalized to higher dimensions, too. It can also
be generalized to any number $m\ge 4$ of types of points such that the
difference between the smallest and the largest translational part
occurring in $\Phi$ is $m^2_{}-1$ (and we really have an LSS). In this case,
it can be shown that $R$ must be at least $m^2-1$, if the corresponding 
substitution system is nicely growing. Anyway, the substitution above
on 4 letters yields the same sequences as the admissible substitution 
\begin{equation*}
a \to dcba, \; b \to dbca, \; c \to cccb, \; d \to dbba, 
\end{equation*}  
showing a modular coincidence in $\Phi_{}^2$. 

\subsection{The Chair Tiling and the Table Tiling}\label{subsec:chairtable}

\begin{figure}[t] 
\begin{tabular}{@{\hspace*{2.55em}}c@{\hspace*{2.5em}}|@{\hspace*{2.5em}}c}
\setlength{\unitlength}{.5cm}
\begin{picture}(10,4)
\put(.5,1){\line(1,0){2}}
\put(.5,1){\line(0,1){2}}
\put(2.5,2){\line(-1,0){1}}
\put(2.5,2){\line(0,-1){1}}
\put(1.5,3){\line(-1,0){1}}
\put(1.5,3){\line(0,-1){1}}

{
\thicklines
\put(3.5,2){\vector(1,0){1.5}}
}

\put(5.5,0){\line(1,0){4}}
\put(5.5,0){\line(0,1){4}}
\put(9.5,2){\line(-1,0){2}}
\put(9.5,2){\line(0,-1){2}}
\put(7.5,4){\line(-1,0){2}}
\put(7.5,4){\line(0,-1){2}}
\put(6.5,1){\line(1,0){2}}
\put(6.5,1){\line(0,1){2}}
\put(8.5,2){\line(0,-1){1}}
\put(7.5,1){\line(0,-1){1}}
\put(5.5,2){\line(1,0){1}}
\put(6.5,3){\line(1,0){1}}
\end{picture}
&
\setlength{\unitlength}{.5cm}
\begin{picture}(10,4)
\put(.5,1.5){\line(1,0){2}}
\put(.5,1.5){\line(0,1){1}}
\put(2.5,2.5){\line(-1,0){2}}
\put(2.5,2.5){\line(0,-1){1}}

{
\thicklines
\put(3,2){\vector(1,0){1.5}}
}

\put(5,1){\line(1,0){4}}
\put(5,1){\line(0,1){2}}
\put(9,3){\line(-1,0){4}}
\put(9,3){\line(0,-1){2}}
\put(6,1){\line(0,1){2}}
\put(8,1){\line(0,1){2}}
\put(6,2){\line(1,0){2}}
\end{picture}
\\
\end{tabular}
\caption{Tile substitution for the Chair Tiling (left) and for the 
Table Tiling (right). \label{fig:chairtable}}
\end{figure}
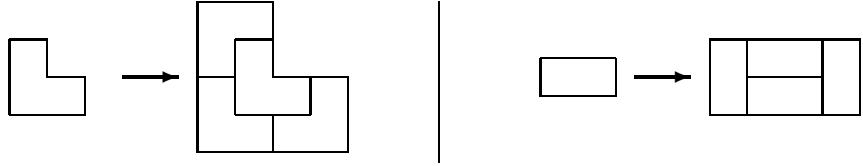

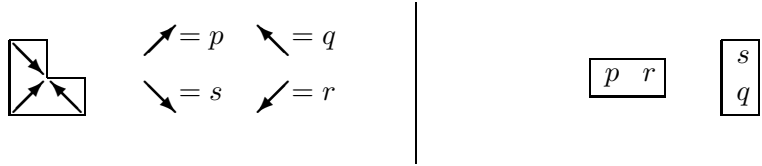
\begin{figure}[t] 
\begin{tabular}{@{\hspace*{0em}}c@{\hspace*{3em}}|@{\hspace*{5em}}c}
\setlength{\unitlength}{.5cm}
\begin{picture}(9,4)
\put(.5,1){\line(1,0){2}}
\put(.5,1){\line(0,1){2}}
\put(2.5,2){\line(-1,0){1}}
\put(2.5,2){\line(0,-1){1}}
\put(1.5,3){\line(-1,0){1}}
\put(1.5,3){\line(0,-1){1}}

{
\thicklines
\put(.6,1.1){\vector(1,1){.8}}
\put(.6,2.9){\vector(1,-1){.8}}
\put(2.4,1.1){\vector(-1,1){.8}}

\put(4.1,2.6){\vector(1,1){.8}}
\put(7.9,2.6){\vector(-1,1){.8}}
\put(4.1,1.9){\vector(1,-1){.8}}
\put(7.9,1.9){\vector(-1,-1){.8}}
}

\put(5,2.9){$= p$}
\put(8,2.9){$= q$}
\put(5,1.4){$= s$}
\put(8,1.4){$= r$}
\end{picture}
&
\setlength{\unitlength}{.5cm}
\begin{picture}(5.5,4)
\put(.5,1.5){\line(1,0){2}}
\put(.5,1.5){\line(0,1){1}}
\put(2.5,2.5){\line(-1,0){2}}
\put(2.5,2.5){\line(0,-1){1}}

\put(.9,1.9){$p$}
\put(1.9,1.9){$r$}

\put(4,1){\line(1,0){1}}
\put(4,1){\line(0,1){2}}
\put(5,3){\line(-1,0){1}}
\put(5,3){\line(0,-1){2}}

\put(4.4,1.4){$q$}
\put(4.4,2.4){$s$}
\end{picture}
\\
\end{tabular}
\caption{Labelling for Chair Tiling (left) and for the Table Tiling
  (right) \label{fig:ctlabel}}  
\end{figure}

The chair tiling and the table tiling (also called the domino tiling)
and their spectral properties are discussed in detail in~\cite{Sol97},
\cite{Rob99} and~\cite{BMS98}. They are generated by the tile substitutions in
Fig.~\ref{fig:chairtable}. The expansion $Q$ acts as scaling by a 
factor $2$, and the rotated versions of the shown tiles are
substituted analogously. Both tile substitutions can be described by
substitutions over the alphabet $\{p,q,r,s\}$ on $L=\Z^2_{}$, where
one uses the labelling of Fig.~\ref{fig:ctlabel}. 

\begin{figure}[t]
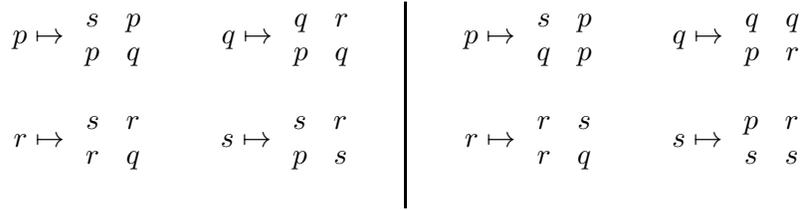
 
\begin{tabular}{@{\hspace*{3.05em}}c@{\hspace*{1.5em}}|@{\hspace*{2em}}c}
$p\mapsto\begin{array}{cc} s & p \\ p & q \end{array}$ \qquad
$q\mapsto\begin{array}{cc} q & r \\ p & q \end{array}$
&
$p\mapsto\begin{array}{cc} s & p \\ q & p \end{array}$ \qquad
$q\mapsto\begin{array}{cc} q & q \\ p & r \end{array}$
\\ \\
$r\mapsto\begin{array}{cc} s & r \\ r & q \end{array}$ \qquad
$s\mapsto\begin{array}{cc} s & r \\ p & s \end{array}$
&
$r\mapsto\begin{array}{cc} r & s \\ r & q \end{array}$ \qquad
$s\mapsto\begin{array}{cc} p & r \\ s & s \end{array}$
\\ \\
\end{tabular}
\caption{Induced substitution for the Chair Tiling (left) 
and for the Table Tiling (right). \label{fig:indsub}}
\end{figure}

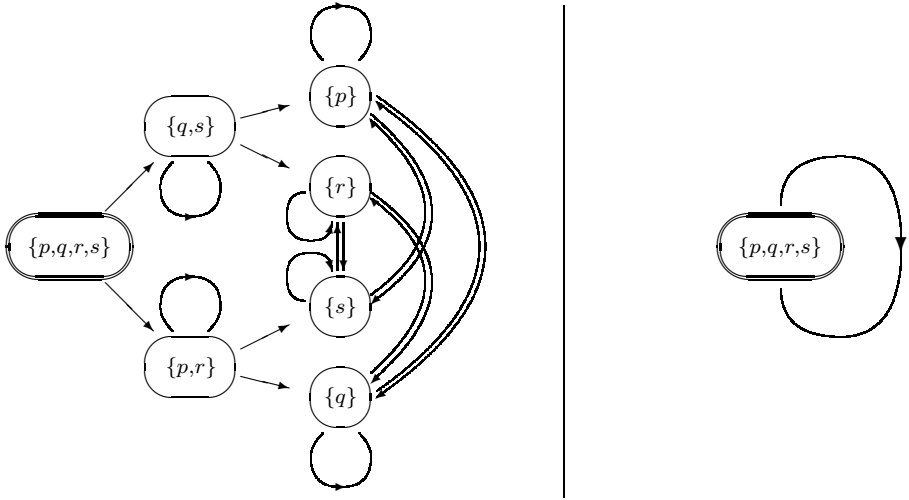
\begin{figure}[t]
\begin{tabular}{@{\hspace*{0em}}c|@{\hspace*{4em}}c@{\hspace*{4.75em}}}
\setlength{\unitlength}{.8cm}
\begin{picture}(9,8)
\put(1,4){\oval(2,1)}
\put(1,4){\oval(2.1,1.1)}
\put(0.3,3.9){$\scriptstyle \{p,q,r,s\}$}

\put(3,6){\oval(1.5,1)}
\put(3,2){\oval(1.5,1)}
\put(2.6,5.9){$\scriptstyle \{q,s\}$}
\put(2.6,1.9){$\scriptstyle \{p,r\}$}

\put(5.5,6.5){\oval(1,1)}
\put(5.5,5){\oval(1,1)}
\put(5.5,3){\oval(1,1)}
\put(5.5,1.5){\oval(1,1)}
\put(5.23,6.4){$\scriptstyle \{p\}$}
\put(5.23,4.9){$\scriptstyle \{r\}$}
\put(5.23,2.9){$\scriptstyle \{s\}$}
\put(5.23,1.4){$\scriptstyle \{q\}$}

\put(1.6,4.6){\vector(1,1){.8}}
\put(1.6,3.4){\vector(1,-1){.8}}

\put(3.85,6.15){\vector(4,1){.8}}
\put(3.85,5.7){\vector(2,-1){.8}}
\put(3.85,1.85){\vector(4,-1){.8}}
\put(3.85,2.3){\vector(2,1){.8}}

\qbezier(2.7,5.4)(2.5,5.2)(2.5,5)
\qbezier(3,4.5)(2.5,4.5)(2.5,5)
\qbezier(3,4.5)(3.5,4.5)(3.5,5)
\qbezier(3.3,5.4)(3.5,5.2)(3.5,5)
\put(3.1,4.5){\vector(1,0){0}}

\qbezier(2.7,2.6)(2.5,2.8)(2.5,3)
\qbezier(3,3.5)(2.5,3.5)(2.5,3)
\qbezier(3,3.5)(3.5,3.5)(3.5,3)
\qbezier(3.3,2.6)(3.5,2.8)(3.5,3)
\put(3.1,3.5){\vector(1,0){0}}

\put(5.45,3.6){\vector(0,1){.8}}
\put(5.55,4.4){\vector(0,-1){.8}}

\qbezier(6.1,6.4)(9.5,4)(6.1,1.6)
\qbezier(6.1,6.5)(9.7,4)(6.1,1.5)
\qbezier(6.0,6.2)(8,4.5)(6.05,3.1)
\qbezier(6.0,6.1)(7.8,4.5)(6.0,3.2)
\qbezier(6.05,1.8)(8,3.5)(6.0,4.9)
\qbezier(6.0,1.9)(7.8,3.5)(6.0,4.8)
\put(6.08,6.43){\vector(-1,1){0}}
\put(6.08,1.47){\vector(-4,-3){0}}
\put(6.0,3.03){\vector(-1,-1){0}}
\put(5.98,6.15){\vector(-2,3){0}}
\put(6.0,1.72){\vector(-1,-1){0}}
\put(5.98,4.83){\vector(-1,1){0}}

\qbezier(4.9,3.1)(4.6,3.1)(4.6,3.5)
\qbezier(4.6,3.5)(4.6,3.9)(5,3.9)
\qbezier(5,3.9)(5.35,3.9)(5.35,3.6)
\put(5.36,3.6){\vector(1,-4){0}}

\qbezier(4.9,4.9)(4.6,4.9)(4.6,4.5)
\qbezier(4.6,4.5)(4.6,4.1)(5,4.1)
\qbezier(5,4.1)(5.35,4.1)(5.35,4.4)
\put(5.36,4.4){\vector(1,4){0}}

\qbezier(5.2,7.1)(5,7.3)(5,7.5)
\qbezier(5.5,8)(5,8)(5,7.5)
\qbezier(5.5,8)(6,8)(6,7.5)
\qbezier(5.8,7.1)(6,7.3)(6,7.5)
\put(5.6,8){\vector(1,0){0}}

\qbezier(5.2,.9)(5,.7)(5,.5)
\qbezier(5.5,0)(5,0)(5,.5)
\qbezier(5.5,0)(6,0)(6,.5)
\qbezier(5.8,.9)(6,.7)(6,.5)
\put(5.6,0){\vector(1,0){0}}
\end{picture}
&
\setlength{\unitlength}{.8cm}
\begin{picture}(4,8)
\put(1.5,4){\oval(2,1)}
\put(1.5,4){\oval(2.1,1.1)}
\put(.8,3.9){$\scriptstyle \{p,q,r,s\}$}

\qbezier(1.5,4.7)(1.5,5.5)(2.5,5.5)
\qbezier(1.5,3.3)(1.5,2.5)(2.5,2.5)
\qbezier(2.5,5.5)(3.5,5.5)(3.5,4)
\qbezier(2.5,2.5)(3.5,2.5)(3.5,4)

{
\thicklines
\put(3.5,3.9){\vector(0,-1){0}}
}
\end{picture}
\\
\end{tabular}
\caption{Coincidence graph for the Chair Tiling (left) and  
for the Table Tiling (right) \label{fig:ctgraph}}
\end{figure}

These labellings yield the (induced) substitutions in
Fig.~\ref{fig:indsub}. Therefore, we obtain the coincidence
graphs shown in Fig.~\ref{fig:ctgraph} (observing that $L=L'=\Z^2_{}$).
Note that the substitution for the table tiling admits no pairwise
coincidence, therefore Corollary \ref{cor:paircoinc} yields that the 
induced LSS does not consist of model sets. 
The substitution is also a bijective substitution, so Corollary
\ref{cor:bij} applies, too. 

These graphs show that the chair tiling is pure point diffractive
(compare to~\cite[Theorem 8.1]{Rob99}, also see~\cite[Example 7.1]{Sol97}),
while the table tiling is not (compare to~\cite[Theorem 8.2]{Rob99}, also
see~\cite[Example 7.3]{Sol97}). We note that the internal space
for the chair tiling is $\Z^2_2$ (which lines up with the inflation
by a factor of $2$, i.e., $\Z^2_2=\lim_{\gets k} \Z^2/2^k_{}\cdot\Z^2_{}$).

\subsection{R1--R2 Tilings}\label{subsec:natalies}
In \cite{PF03}, a family of substitutions is considered which fulfil
two conditions, namely, (R1) and (R2). In our notation, these read as
follows. Let $(\V,\Phi)$ be an admissible primitive LSS
supported on $L=\Z^d$, $\V=(V_1,\ldots,V_{2m})$, $[L:L']=m$.
\begin{itemize}
\item[(R1)] $\Psi^{}_0[a_i] = \{i, i+m\}$ ($1 \le i \le m$), where $a_i$
runs through $L/L'$. 
\item[(R2)] $\Phi_{j,i} \cap \Phi_{j,i+m} = \varnothing$ for all $i,j$. 
\end{itemize}
The first condition tells us that no $\Psi^{}_0[a]$ is a singleton
set: each one contains two elements $i, i+m$. Therefore, the second
condition tells us that there is no pairwise coincidence 
in any $\Psi^{}_0[a]$. Thus, it follows from Corollary $\ref{cor:paircoinc}$ 
that none of the multi-component Delone sets (or tilings) induced by these
substitutions consist of model sets.

\subsection{The Semi-Detached House Tiling}\label{subsec:haus}
\begin{figure}[t] 
\epsfig{file=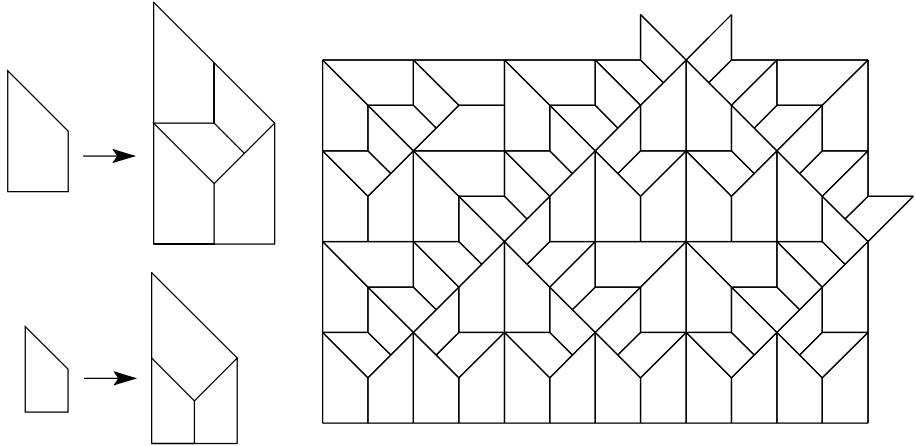}
\caption{The substitution for the Semi-Detached House Tiling (left)
  and a part of the tiling, which (hopefully) illustrates its
  name. \label{fig:haussub}}  
\end{figure}

The tile substitution in Fig.~\ref{fig:haussub} (with expansion
factor $2$) appears as unsolved problem 
in~\cite[Abbildung 17]{Fre02}. As above, rotated and reflected versions
of the tiles are substituted analogously to the shown tiles.

\begin{figure}[t] 
\begin{center}
\setlength{\unitlength}{1cm}
\begin{picture}(14,6)
{
\thicklines

\put(0,3.5){\line(1,1){1}}
\put(0,3.5){\line(0,1){2}}
\put(1,5.5){\line(0,-1){1}}
\put(1,5.5){\line(-1,0){1}}

\put(1,.5){\line(-1,0){1}}
\put(1,.5){\line(0,1){2}}
\put(0,1.5){\line(0,-1){1}}
\put(0,1.5){\line(1,1){1}}

\put(2.5,3.5){\line(-1,1){1}}
\put(2.5,3.5){\line(0,1){2}}
\put(1.5,5.5){\line(0,-1){1}}
\put(1.5,5.5){\line(1,0){1}}

\put(1.5,.5){\line(1,0){1}}
\put(1.5,.5){\line(0,1){2}}
\put(2.5,1.5){\line(0,-1){1}}
\put(2.5,1.5){\line(-1,1){1}}

\put(3,0){\line(1,0){2}}
\put(3,0){\line(0,1){1}}
\put(4,1){\line(-1,0){1}}
\put(4,1){\line(1,-1){1}}

\put(3,2.5){\line(1,0){2}}
\put(3,2.5){\line(0,-1){1}}
\put(4,1.5){\line(-1,0){1}}
\put(4,1.5){\line(1,1){1}}

\put(5,3.5){\line(-1,0){2}}
\put(5,3.5){\line(0,1){1}}
\put(4,4.5){\line(1,0){1}}
\put(4,4.5){\line(-1,-1){1}}

\put(5,6){\line(-1,0){2}}
\put(5,6){\line(0,-1){1}}
\put(4,5){\line(1,0){1}}
\put(4,5){\line(-1,1){1}}

\put(6,4){\line(0,1){2}}
\put(6,4){\line(1,-1){1}}
\put(8,4){\line(-1,-1){1}}
\put(8,4){\line(-1,1){2}}

\put(6.5,6){\line(1,-1){2}}
\put(6.5,6){\line(1,0){2}}
\put(9.5,5){\line(-1,1){1}}
\put(9.5,5){\line(-1,-1){1}}

\put(6,1){\line(1,1){1}}
\put(6,1){\line(1,-1){1}}
\put(9,0){\line(-1,1){2}}
\put(9,0){\line(-1,0){2}}

\put(9.5,0){\line(-1,1){2}}
\put(9.5,0){\line(0,1){2}}
\put(8.5,3){\line(-1,-1){1}}
\put(8.5,3){\line(1,-1){1}}

\put(10.5,5){\line(1,1){1}}
\put(10.5,5){\line(1,-1){1}}
\put(13.5,6){\line(-1,0){2}}
\put(13.5,6){\line(-1,-1){2}}

\put(10.5,0){\line(1,1){2}}
\put(10.5,0){\line(0,1){2}}
\put(11.5,3){\line(-1,-1){1}}
\put(11.5,3){\line(1,-1){1}}

\put(14,6){\line(0,-1){2}}
\put(14,6){\line(-1,-1){2}}
\put(13,3){\line(1,1){1}}
\put(13,3){\line(-1,1){1}}

\put(11,0){\line(1,1){2}}
\put(11,0){\line(1,0){2}}
\put(14,1){\line(-1,1){1}}
\put(14,1){\line(-1,-1){1}}
}

\put(1,4.5){\line(-1,1){1}}
\put(1,4.5){\line(-1,0){1}}

\put(0,1.5){\line(1,-1){1}}
\put(0,1.5){\line(1,0){1}}

\put(1.5,4.5){\line(1,1){1}}
\put(1.5,4.5){\line(1,0){1}}

\put(2.5,1.5){\line(-1,-1){1}}
\put(2.5,1.5){\line(-1,0){1}}

\put(4,1){\line(-1,-1){1}}
\put(4,1){\line(0,-1){1}}

\put(4,1.5){\line(-1,1){1}}
\put(4,1.5){\line(0,1){1}}

\put(4,4.5){\line(1,-1){1}}
\put(4,4.5){\line(0,-1){1}}

\put(4,5){\line(0,1){1}}
\put(4,5){\line(1,1){1}}

\put(6,4){\line(1,0){2}}
\put(6,4){\line(1,1){1}}
\put(6,5){\line(1,0){1}}
\put(7,3){\line(0,1){2}}

\put(7.5,5){\line(1,0){2}}
\put(7.5,5){\line(0,1){1}}
\put(8.5,6){\line(0,-1){2}}
\put(8.5,6){\line(-1,-1){1}}

\put(8,1){\line(-1,0){2}}
\put(8,1){\line(0,-1){1}}
\put(8,1){\line(-1,-1){1}}
\put(7,0){\line(0,1){2}}

\put(8.5,1){\line(0,1){2}}
\put(8.5,1){\line(1,0){1}}
\put(8.5,1){\line(1,1){1}}
\put(7.5,2){\line(1,0){2}}

\put(12.5,5){\line(-1,0){2}}
\put(12.5,5){\line(0,1){1}}
\put(12.5,5){\line(-1,1){1}}
\put(11.5,4){\line(0,1){2}}

\put(13,5){\line(1,0){1}}
\put(13,5){\line(1,-1){1}}
\put(13,5){\line(0,-1){2}}
\put(12,4){\line(1,0){2}}

\put(11.5,1){\line(-1,0){1}}
\put(11.5,1){\line(0,1){2}}
\put(11.5,1){\line(-1,1){1}}
\put(10.5,2){\line(1,0){2}}

\put(12,1){\line(1,0){2}}
\put(12,1){\line(0,-1){1}}
\put(12,1){\line(1,-1){1}}
\put(13,0){\line(0,1){2}}

\put(2.0,.7){$A$}
\put(1.6,1.1){$A$}
\put(1.6,1.6){$A$}

\put(.15,.65){$\overline{A}$}
\put(.55,1.05){$\overline{A}$}
\put(.55,1.6){$\overline{A}$}

\put(3.6,5.6){$B$}
\put(4.2,5.6){$B$}
\put(4.6,5.1){$B$}

\put(3.6,3.6){$\overline{B}$}
\put(4.2,3.6){$\overline{B}$}
\put(4.6,4.05){$\overline{B}$}

\put(3.2,0.6){$C$}
\put(3.6,0.1){$C$}
\put(4.2,0.1){$C$}

\put(3.2,1.6){$\overline{C}$}
\put(3.6,2.05){$\overline{C}$}
\put(4.2,2.05){$\overline{C}$}

\put(1.6,5.0){$D$}
\put(2.0,4.6){$D$}
\put(2.0,4.05){$D$}

\put(.55,5.0){$\overline{D}$}
\put(.15,4.6){$\overline{D}$}
\put(.15,4.0){$\overline{D}$}

\put(8.6,1.6){$V$}
\put(8.1,1.6){$V$}
\put(8.6,2.2){$V$}
\put(8.1,2.2){$V$}
\put(9.1,1.15){$V$}
\put(9.1,0.55){$V$}

\put(11.6,1.55){$\overline{V}$}
\put(11.1,1.55){$\overline{V}$}
\put(11.6,2.15){$\overline{V}$}
\put(11.1,2.15){$\overline{V}$}
\put(10.6,1.15){$\overline{V}$}
\put(10.6,0.55){$\overline{V}$}

\put(7.1,0.6){$X$}
\put(6.6,0.6){$X$}
\put(7.1,1.2){$X$}
\put(6.6,1.2){$X$}
\put(7.6,0.15){$X$}
\put(8.1,0.15){$X$}

\put(11.6,4.55){$\overline{X}$}
\put(11.1,4.55){$\overline{X}$}
\put(11.6,5.15){$\overline{X}$}
\put(11.1,5.15){$\overline{X}$}
\put(12.1,5.55){$\overline{X}$}
\put(12.6,5.55){$\overline{X}$}

\put(8.6,4.6){$Y$}
\put(8.1,4.6){$Y$}
\put(8.6,5.2){$Y$}
\put(8.1,5.2){$Y$}
\put(7.6,5.6){$Y$}
\put(7.1,5.6){$Y$}

\put(13.1,0.55){$\overline{Y}$}
\put(12.6,0.55){$\overline{Y}$}
\put(13.1,1.15){$\overline{Y}$}
\put(12.6,1.15){$\overline{Y}$}
\put(12.1,0.1){$\overline{Y}$}
\put(11.6,0.1){$\overline{Y}$}

\put(7.1,3.6){$Z$}
\put(6.6,3.6){$Z$}
\put(7.1,4.2){$Z$}
\put(6.6,4.2){$Z$}
\put(6.1,4.6){$Z$}
\put(6.1,5.1){$Z$}

\put(13.1,3.55){$\overline{Z}$}
\put(12.6,3.55){$\overline{Z}$}
\put(13.1,4.15){$\overline{Z}$}
\put(12.6,4.15){$\overline{Z}$}
\put(13.6,4.55){$\overline{Z}$}
\put(13.6,5.1){$\overline{Z}$}
\end{picture}
\end{center}
\caption{Labelling of the Semi-Detached House Tiling
 \label{fig:hauslabel}}
\end{figure}
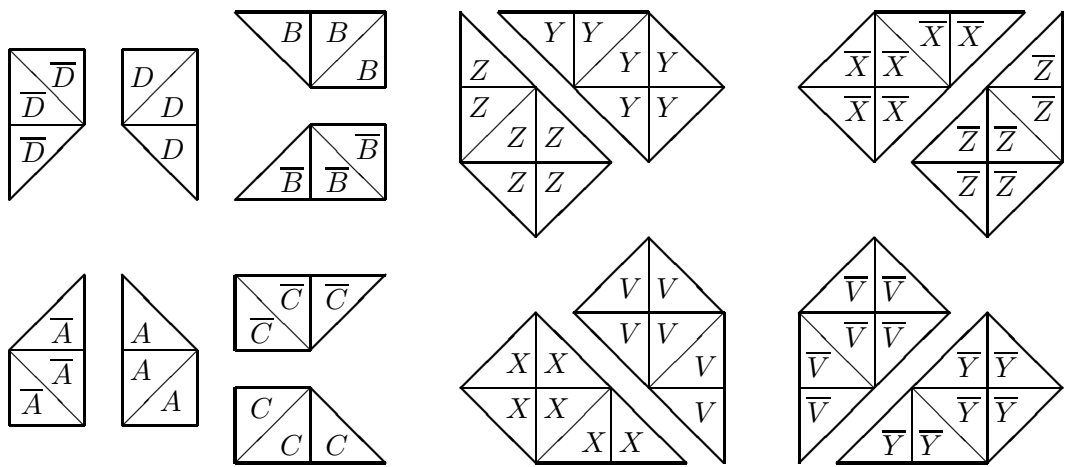

Each of the two semi-houses (i.e., the big and the small one, see
tile substitution in Fig.~\ref{fig:haussub}) arises in $4$ different
orientations and its reflected versions. To derive a substitution, we
use the labelling in Fig.~\ref{fig:hauslabel}. 

\begin{figure}[t] 
\begin{tabular}{c@{\qquad}c}
\setlength{\unitlength}{.65cm}
\begin{picture}(10.5,25)
\put(0,23.8){$a$)}
\put(1,23.5){\line(1,0){1}}
\put(1,23.5){\line(1,1){1}}
\put(1,23.5){\line(0,1){1}}
\put(2,24.5){\line(-1,0){1}}
\put(2,24.5){\line(0,-1){1}}
\put(1.6,23.6){$\scriptstyle A$}
\put(1.1,24.05){$\scriptstyle A$}
\put(2.15,23.55){,}
\put(2.5,23.5){\line(1,0){1}}
\put(2.5,23.5){\line(1,1){1}}
\put(2.5,23.5){\line(0,1){1}}
\put(3.5,24.5){\line(-1,0){1}}
\put(3.5,24.5){\line(0,-1){1}}
\put(3.1,23.6){$\scriptstyle Z$}
\put(2.6,24.05){$\scriptstyle Z$}
{
\thicklines
\put(4,24){\vector(1,0){1}}
}
\put(5.5,23){\line(1,0){2}}
\put(5.5,23){\line(1,1){2}}
\put(5.5,23){\line(0,1){2}}
\put(7.5,25){\line(-1,0){2}}
\put(7.5,25){\line(0,-1){2}}
\put(5.5,25){\line(1,-1){2}}
\put(6.5,23){\line(0,1){2}}
\put(5.5,24){\line(1,0){2}}
\put(7.1,24.1){$\scriptstyle \overline{A}$}
\put(6.6,24.55){$\scriptstyle Z$}
\put(6.1,24.55){$\scriptstyle Z$}
\put(5.6,24.1){$\scriptstyle A$}
\put(5.6,23.55){$\scriptstyle A$}
\put(6.1,23.1){$\scriptstyle A$}
\put(6.6,23.1){$\scriptstyle \overline{A}$}
\put(7.1,23.55){$\scriptstyle \overline{A}$}

\put(0,21.5){$b$)}
\put(1,21.0){\line(1,0){1}}
\put(1,21.0){\line(1,1){1}}
\put(1,21.0){\line(0,1){1}}
\put(2,22.0){\line(-1,0){1}}
\put(2,22.0){\line(0,-1){1}}
\put(1.6,21.1){$\scriptstyle B$}
\put(1.1,21.55){$\scriptstyle B$}
\put(2.15,21.05){,}
\put(2.5,21.0){\line(1,0){1}}
\put(2.5,21.0){\line(1,1){1}}
\put(2.5,21.0){\line(0,1){1}}
\put(3.5,22.0){\line(-1,0){1}}
\put(3.5,22.0){\line(0,-1){1}}
\put(3.1,21.1){$\scriptstyle Y$}
\put(2.6,21.55){$\scriptstyle Y$}
{
\thicklines
\put(4,21.5){\vector(1,0){1}}
}
\put(5.5,20.5){\line(1,0){2}}
\put(5.5,20.5){\line(1,1){2}}
\put(5.5,20.5){\line(0,1){2}}
\put(7.5,22.5){\line(-1,0){2}}
\put(7.5,22.5){\line(0,-1){2}}
\put(5.5,22.5){\line(1,-1){2}}
\put(6.5,20.5){\line(0,1){2}}
\put(5.5,21.5){\line(1,0){2}}
\put(7.1,21.6){$\scriptstyle B$}
\put(6.6,22.05){$\scriptstyle B$}
\put(6.1,22.05){$\scriptstyle B$}
\put(5.6,21.6){$\scriptstyle Y$}
\put(5.6,21.05){$\scriptstyle Y$}
\put(6.1,20.6){$\scriptstyle \overline{B}$}
\put(6.6,20.6){$\scriptstyle \overline{B}$}
\put(7.1,21.05){$\scriptstyle \overline{B}$}

\put(0,18.8){$c$)}
\put(1,18.5){\line(1,0){1}}
\put(1,18.5){\line(1,1){1}}
\put(1,18.5){\line(0,1){1}}
\put(2,19.5){\line(-1,0){1}}
\put(2,19.5){\line(0,-1){1}}
\put(1.6,18.6){$\scriptstyle C$}
\put(1.1,19.05){$\scriptstyle C$}
\put(2.15,18.55){,}
\put(2.5,18.5){\line(1,0){1}}
\put(2.5,18.5){\line(1,1){1}}
\put(2.5,18.5){\line(0,1){1}}
\put(3.5,19.5){\line(-1,0){1}}
\put(3.5,19.5){\line(0,-1){1}}
\put(3.1,18.6){$\scriptstyle X$}
\put(2.6,19.05){$\scriptstyle X$}
{
\thicklines
\put(4,19){\vector(1,0){1}}
}
\put(5.5,18){\line(1,0){2}}
\put(5.5,18){\line(1,1){2}}
\put(5.5,18){\line(0,1){2}}
\put(7.5,20){\line(-1,0){2}}
\put(7.5,20){\line(0,-1){2}}
\put(5.5,20){\line(1,-1){2}}
\put(6.5,18){\line(0,1){2}}
\put(5.5,19){\line(1,0){2}}
\put(7.1,19.1){$\scriptstyle X$}
\put(6.6,19.55){$\scriptstyle \overline{C}$}
\put(6.1,19.55){$\scriptstyle \overline{C}$}
\put(5.6,19.1){$\scriptstyle \overline{C}$}
\put(5.6,18.55){$\scriptstyle C$}
\put(6.1,18.1){$\scriptstyle C$}
\put(6.6,18.1){$\scriptstyle C$}
\put(7.1,18.55){$\scriptstyle X$}

\put(0,16.5){$d$)}
\put(1,16.0){\line(1,0){1}}
\put(1,16.0){\line(1,1){1}}
\put(1,16.0){\line(0,1){1}}
\put(2,17.0){\line(-1,0){1}}
\put(2,17.0){\line(0,-1){1}}
\put(1.6,16.1){$\scriptstyle D$}
\put(1.1,16.55){$\scriptstyle D$}
\put(2.15,16.05){,}
\put(2.5,16.0){\line(1,0){1}}
\put(2.5,16.0){\line(1,1){1}}
\put(2.5,16.0){\line(0,1){1}}
\put(3.5,17.0){\line(-1,0){1}}
\put(3.5,17.0){\line(0,-1){1}}
\put(3.1,16.1){$\scriptstyle V$}
\put(2.6,16.55){$\scriptstyle V$}
{
\thicklines
\put(4,16.5){\vector(1,0){1}}
}
\put(5.5,15.5){\line(1,0){2}}
\put(5.5,15.5){\line(1,1){2}}
\put(5.5,15.5){\line(0,1){2}}
\put(7.5,17.5){\line(-1,0){2}}
\put(7.5,17.5){\line(0,-1){2}}
\put(5.5,17.5){\line(1,-1){2}}
\put(6.5,15.5){\line(0,1){2}}
\put(5.5,16.5){\line(1,0){2}}
\put(7.1,16.6){$\scriptstyle D$}
\put(6.6,17.05){$\scriptstyle D$}
\put(6.1,17.05){$\scriptstyle \overline{D}$}
\put(5.6,16.6){$\scriptstyle \overline{D}$}
\put(5.6,16.05){$\scriptstyle \overline{D}$}
\put(6.1,15.6){$\scriptstyle V$}
\put(6.6,15.6){$\scriptstyle V$}
\put(7.1,16.05){$\scriptstyle D$}

\put(0,13.8){$e$)}
\put(1,13.5){\line(1,0){1}}
\put(1,13.5){\line(1,1){1}}
\put(1,13.5){\line(0,1){1}}
\put(2,14.5){\line(-1,0){1}}
\put(2,14.5){\line(0,-1){1}}
\put(1.6,13.6){$\scriptstyle \overline{A}$}
\put(1.1,14.05){$\scriptstyle Z$}
\put(2.15,13.55){,}
\put(2.5,13.5){\line(1,0){1}}
\put(2.5,13.5){\line(1,1){1}}
\put(2.5,13.5){\line(0,1){1}}
\put(3.5,14.5){\line(-1,0){1}}
\put(3.5,14.5){\line(0,-1){1}}
\put(3.1,13.6){$\scriptstyle \overline{Z}$}
\put(2.6,14.05){$\scriptstyle Z$}
\put(3.65,13.55){,}
\put(4,13.5){\line(1,0){1}}
\put(4,13.5){\line(1,1){1}}
\put(4,13.5){\line(0,1){1}}
\put(5,14.5){\line(-1,0){1}}
\put(5,14.5){\line(0,-1){1}}
\put(4.6,13.6){$\scriptstyle \overline{Z}$}
\put(4.1,14.05){$\scriptstyle Y$}
{
\thicklines
\put(5.5,14){\vector(1,0){1}}
}
\put(7.0,13){\line(1,0){2}}
\put(7.0,13){\line(1,1){2}}
\put(7.0,13){\line(0,1){2}}
\put(9.0,15){\line(-1,0){2}}
\put(9.0,15){\line(0,-1){2}}
\put(7.0,15){\line(1,-1){2}}
\put(8.0,13){\line(0,1){2}}
\put(7.0,14){\line(1,0){2}}
\put(8.6,14.1){$\scriptstyle \overline{Z}$}
\put(8.1,14.55){$\scriptstyle Z$}
\put(7.6,14.55){$\scriptstyle Z$}
\put(7.1,14.1){$\scriptstyle Y$}
\put(7.1,13.55){$\scriptstyle Y$}
\put(7.6,13.1){$\scriptstyle \overline{Z}$}
\put(8.1,13.1){$\scriptstyle \overline{Z}$}
\put(8.6,13.55){$\scriptstyle \overline{Z}$}

\put(0,11.3){$f$)}
\put(1,11){\line(1,0){1}}
\put(1,11){\line(1,1){1}}
\put(1,11){\line(0,1){1}}
\put(2,12){\line(-1,0){1}}
\put(2,12){\line(0,-1){1}}
\put(1.6,11.1){$\scriptstyle \overline{B}$}
\put(1.1,11.55){$\scriptstyle Y$}
\put(2.15,11.05){,}
\put(2.5,11){\line(1,0){1}}
\put(2.5,11){\line(1,1){1}}
\put(2.5,11){\line(0,1){1}}
\put(3.5,12){\line(-1,0){1}}
\put(3.5,12){\line(0,-1){1}}
\put(3.1,11.1){$\scriptstyle \overline{Y}$}
\put(2.6,11.55){$\scriptstyle Y$}
\put(3.65,11.05){,}
\put(4,11){\line(1,0){1}}
\put(4,11){\line(1,1){1}}
\put(4,11){\line(0,1){1}}
\put(5,12){\line(-1,0){1}}
\put(5,12){\line(0,-1){1}}
\put(4.6,11.1){$\scriptstyle \overline{Y}$}
\put(4.1,11.55){$\scriptstyle Z$}
{
\thicklines
\put(5.5,11.5){\vector(1,0){1}}
}
\put(7.0,10.5){\line(1,0){2}}
\put(7.0,10.5){\line(1,1){2}}
\put(7.0,10.5){\line(0,1){2}}
\put(9.0,12.5){\line(-1,0){2}}
\put(9.0,12.5){\line(0,-1){2}}
\put(7.0,12.5){\line(1,-1){2}}
\put(8.0,10.5){\line(0,1){2}}
\put(7.0,11.5){\line(1,0){2}}
\put(8.6,11.6){$\scriptstyle \overline{Y}$}
\put(8.1,12.05){$\scriptstyle Z$}
\put(7.6,12.05){$\scriptstyle Z$}
\put(7.1,11.6){$\scriptstyle Y$}
\put(7.1,11.05){$\scriptstyle Y$}
\put(7.6,10.6){$\scriptstyle \overline{Y}$}
\put(8.1,10.6){$\scriptstyle \overline{Y}$}
\put(8.6,11.05){$\scriptstyle \overline{Y}$}

\put(0,8.8){$g$)}
\put(1,8.5){\line(1,0){1}}
\put(1,8.5){\line(1,1){1}}
\put(1,8.5){\line(0,1){1}}
\put(2,9.5){\line(-1,0){1}}
\put(2,9.5){\line(0,-1){1}}
\put(1.6,8.6){$\scriptstyle X$}
\put(1.1,9.05){$\scriptstyle \overline{C}$}
\put(2.15,8.55){,}
\put(2.5,8.5){\line(1,0){1}}
\put(2.5,8.5){\line(1,1){1}}
\put(2.5,8.5){\line(0,1){1}}
\put(3.5,9.5){\line(-1,0){1}}
\put(3.5,9.5){\line(0,-1){1}}
\put(3.1,8.6){$\scriptstyle X$}
\put(2.6,9.05){$\scriptstyle \overline{X}$}
\put(3.65,8.55){,}
\put(4,8.5){\line(1,0){1}}
\put(4,8.5){\line(1,1){1}}
\put(4,8.5){\line(0,1){1}}
\put(5,9.5){\line(-1,0){1}}
\put(5,9.5){\line(0,-1){1}}
\put(4.6,8.6){$\scriptstyle V$}
\put(4.1,9.05){$\scriptstyle \overline{X}$}
{
\thicklines
\put(5.5,9){\vector(1,0){1}}
}
\put(7.0,8){\line(1,0){2}}
\put(7.0,8){\line(1,1){2}}
\put(7.0,8){\line(0,1){2}}
\put(9.0,10){\line(-1,0){2}}
\put(9.0,10){\line(0,-1){2}}
\put(7.0,10){\line(1,-1){2}}
\put(8.0,8){\line(0,1){2}}
\put(7.0,9){\line(1,0){2}}
\put(8.6,9.1){$\scriptstyle X$}
\put(8.1,9.55){$\scriptstyle \overline{X}$}
\put(7.6,9.55){$\scriptstyle \overline{X}$}
\put(7.1,9.1){$\scriptstyle \overline{X}$}
\put(7.1,8.55){$\scriptstyle \overline{X}$}
\put(7.6,8.1){$\scriptstyle V$}
\put(8.1,8.1){$\scriptstyle V$}
\put(8.6,8.55){$\scriptstyle X$}

\put(0,6.3){$h$)}
\put(1,6){\line(1,0){1}}
\put(1,6){\line(1,1){1}}
\put(1,6){\line(0,1){1}}
\put(2,7){\line(-1,0){1}}
\put(2,7){\line(0,-1){1}}
\put(1.6,6.1){$\scriptstyle V$}
\put(1.1,6.55){$\scriptstyle \overline{D}$}
\put(2.15,6.05){,}
\put(2.5,6){\line(1,0){1}}
\put(2.5,6){\line(1,1){1}}
\put(2.5,6){\line(0,1){1}}
\put(3.5,7){\line(-1,0){1}}
\put(3.5,7){\line(0,-1){1}}
\put(3.1,6.1){$\scriptstyle V$}
\put(2.6,6.55){$\scriptstyle \overline{V}$}
\put(3.65,6.05){,}
\put(4,6){\line(1,0){1}}
\put(4,6){\line(1,1){1}}
\put(4,6){\line(0,1){1}}
\put(5,7){\line(-1,0){1}}
\put(5,7){\line(0,-1){1}}
\put(4.6,6.1){$\scriptstyle X$}
\put(4.1,6.55){$\scriptstyle \overline{V}$}
{
\thicklines
\put(5.5,6.5){\vector(1,0){1}}
}
\put(7.0,5.5){\line(1,0){2}}
\put(7.0,5.5){\line(1,1){2}}
\put(7.0,5.5){\line(0,1){2}}
\put(9.0,7.5){\line(-1,0){2}}
\put(9.0,7.5){\line(0,-1){2}}
\put(7.0,7.5){\line(1,-1){2}}
\put(8.0,5.5){\line(0,1){2}}
\put(7.0,6.5){\line(1,0){2}}
\put(8.6,6.6){$\scriptstyle X$}
\put(8.1,7.05){$\scriptstyle \overline{V}$}
\put(7.6,7.05){$\scriptstyle \overline{V}$}
\put(7.1,6.6){$\scriptstyle \overline{V}$}
\put(7.1,6.05){$\scriptstyle \overline{V}$}
\put(7.6,5.6){$\scriptstyle V$}
\put(8.1,5.6){$\scriptstyle V$}
\put(8.6,6.05){$\scriptstyle X$}

\put(0,3.8){$i$)}
\put(1,3.5){\line(1,0){1}}
\put(1,3.5){\line(1,1){1}}
\put(1,3.5){\line(0,1){1}}
\put(2,4.5){\line(-1,0){1}}
\put(2,4.5){\line(0,-1){1}}
\put(1.6,3.6){$\scriptstyle \overline{Z}$}
\put(1.1,4.05){$\scriptstyle \overline{X}$}
\put(2.15,3.55){,}
\put(2.5,3.5){\line(1,0){1}}
\put(2.5,3.5){\line(1,1){1}}
\put(2.5,3.5){\line(0,1){1}}
\put(3.5,4.5){\line(-1,0){1}}
\put(3.5,4.5){\line(0,-1){1}}
\put(3.1,3.6){$\scriptstyle \overline{X}$}
\put(2.6,4.05){$\scriptstyle \overline{Z}$}
\put(3.65,3.55){,}
\put(4,3.5){\line(1,0){1}}
\put(4,3.5){\line(1,1){1}}
\put(4,3.5){\line(0,1){1}}
\put(5,4.5){\line(-1,0){1}}
\put(5,4.5){\line(0,-1){1}}
\put(4.6,3.6){$\scriptstyle \overline{A}$}
\put(4.1,4.05){$\scriptstyle \overline{Z}$}
\put(5.15,3.55){,}
\put(5.5,3.5){\line(1,0){1}}
\put(5.5,3.5){\line(1,1){1}}
\put(5.5,3.5){\line(0,1){1}}
\put(6.5,4.5){\line(-1,0){1}}
\put(6.5,4.5){\line(0,-1){1}}
\put(6.1,3.6){$\scriptstyle \overline{X}$}
\put(5.6,4.05){$\scriptstyle \overline{C}$}
{
\thicklines
\put(7.0,4){\vector(1,0){1}}
}
\put(8.5,3){\line(1,0){2}}
\put(8.5,3){\line(1,1){2}}
\put(8.5,3){\line(0,1){2}}
\put(10.5,5){\line(-1,0){2}}
\put(10.5,5){\line(0,-1){2}}
\put(8.5,5){\line(1,-1){2}}
\put(9.5,3){\line(0,1){2}}
\put(8.5,4){\line(1,0){2}}
\put(10.1,4.1){$\scriptstyle \overline{Z}$}
\put(9.6,4.55){$\scriptstyle \overline{X}$}
\put(9.1,4.55){$\scriptstyle \overline{X}$}
\put(8.6,4.1){$\scriptstyle \overline{X}$}
\put(8.6,3.55){$\scriptstyle \overline{X}$}
\put(9.1,3.1){$\scriptstyle \overline{Z}$}
\put(9.6,3.1){$\scriptstyle \overline{Z}$}
\put(10.1,3.55){$\scriptstyle \overline{Z}$}

\put(0,1.3){$k$)}
\put(1,1){\line(1,0){1}}
\put(1,1){\line(1,1){1}}
\put(1,1){\line(0,1){1}}
\put(2,2){\line(-1,0){1}}
\put(2,2){\line(0,-1){1}}
\put(1.6,1.1){$\scriptstyle \overline{Y}$}
\put(1.1,1.55){$\scriptstyle \overline{V}$}
\put(2.15,1.05){,}
\put(2.5,1){\line(1,0){1}}
\put(2.5,1){\line(1,1){1}}
\put(2.5,1){\line(0,1){1}}
\put(3.5,2){\line(-1,0){1}}
\put(3.5,2){\line(0,-1){1}}
\put(3.1,1.1){$\scriptstyle \overline{V}$}
\put(2.6,1.55){$\scriptstyle \overline{Y}$}
\put(3.65,1.05){,}
\put(4,1){\line(1,0){1}}
\put(4,1){\line(1,1){1}}
\put(4,1){\line(0,1){1}}
\put(5,2){\line(-1,0){1}}
\put(5,2){\line(0,-1){1}}
\put(4.6,1.1){$\scriptstyle \overline{B}$}
\put(4.1,1.55){$\scriptstyle \overline{Y}$}
\put(5.15,1.05){,}
\put(5.5,1){\line(1,0){1}}
\put(5.5,1){\line(1,1){1}}
\put(5.5,1){\line(0,1){1}}
\put(6.5,2){\line(-1,0){1}}
\put(6.5,2){\line(0,-1){1}}
\put(6.1,1.1){$\scriptstyle \overline{V}$}
\put(5.6,1.55){$\scriptstyle \overline{D}$}
{
\thicklines
\put(7.0,1.5){\vector(1,0){1}}
}
\put(8.5,0.5){\line(1,0){2}}
\put(8.5,0.5){\line(1,1){2}}
\put(8.5,0.5){\line(0,1){2}}
\put(10.5,2.5){\line(-1,0){2}}
\put(10.5,2.5){\line(0,-1){2}}
\put(8.5,2.5){\line(1,-1){2}}
\put(9.5,0.5){\line(0,1){2}}
\put(8.5,1.5){\line(1,0){2}}
\put(10.1,1.6){$\scriptstyle \overline{Y}$}
\put(9.6,2.05){$\scriptstyle \overline{V}$}
\put(9.1,2.05){$\scriptstyle \overline{V}$}
\put(8.6,1.6){$\scriptstyle \overline{V}$}
\put(8.6,1.05){$\scriptstyle \overline{V}$}
\put(9.1,0.6){$\scriptstyle \overline{Y}$}
\put(9.6,0.6){$\scriptstyle \overline{Y}$}
\put(10.1,1.05){$\scriptstyle \overline{Y}$}
\end{picture}
&
\setlength{\unitlength}{.65cm}
\begin{picture}(10.5,25)
\put(0,23.8){$\overline{a}$)}
\put(1,23.5){\line(1,0){1}}
\put(1,24.5){\line(1,-1){1}}
\put(1,23.5){\line(0,1){1}}
\put(2,24.5){\line(-1,0){1}}
\put(2,24.5){\line(0,-1){1}}
\put(1.1,23.6){$\scriptstyle \overline{A}$}
\put(1.6,24.05){$\scriptstyle \overline{A}$}
\put(2.15,23.55){,}
\put(2.5,23.5){\line(1,0){1}}
\put(2.5,24.5){\line(1,-1){1}}
\put(2.5,23.5){\line(0,1){1}}
\put(3.5,24.5){\line(-1,0){1}}
\put(3.5,24.5){\line(0,-1){1}}
\put(2.6,23.6){$\scriptstyle \overline{Z}$}
\put(3.1,24.05){$\scriptstyle \overline{Z}$}
{
\thicklines
\put(4,24){\vector(1,0){1}}
}
\put(5.5,23){\line(1,0){2}}
\put(5.5,23){\line(1,1){2}}
\put(5.5,23){\line(0,1){2}}
\put(7.5,25){\line(-1,0){2}}
\put(7.5,25){\line(0,-1){2}}
\put(5.5,25){\line(1,-1){2}}
\put(6.5,23){\line(0,1){2}}
\put(5.5,24){\line(1,0){2}}
\put(7.1,24.1){$\scriptstyle \overline{A}$}
\put(6.6,24.55){$\scriptstyle \overline{Z}$}
\put(6.1,24.55){$\scriptstyle \overline{Z}$}
\put(5.6,24.1){$\scriptstyle A$}
\put(5.6,23.55){$\scriptstyle A$}
\put(6.1,23.1){$\scriptstyle A$}
\put(6.6,23.1){$\scriptstyle \overline{A}$}
\put(7.1,23.55){$\scriptstyle \overline{A}$}

\put(0,21.5){$\overline{b}$)}
\put(1,21.0){\line(1,0){1}}
\put(1,22.0){\line(1,-1){1}}
\put(1,21.0){\line(0,1){1}}
\put(2,22.0){\line(-1,0){1}}
\put(2,22.0){\line(0,-1){1}}
\put(1.1,21.1){$\scriptstyle \overline{B}$}
\put(1.6,21.55){$\scriptstyle \overline{B}$}
\put(2.15,21.05){,}
\put(2.5,21.0){\line(1,0){1}}
\put(2.5,22.0){\line(1,-1){1}}
\put(2.5,21.0){\line(0,1){1}}
\put(3.5,22.0){\line(-1,0){1}}
\put(3.5,22.0){\line(0,-1){1}}
\put(2.6,21.1){$\scriptstyle \overline{Y}$}
\put(3.1,21.55){$\scriptstyle \overline{Y}$}
{
\thicklines
\put(4,21.5){\vector(1,0){1}}
}
\put(5.5,20.5){\line(1,0){2}}
\put(5.5,20.5){\line(1,1){2}}
\put(5.5,20.5){\line(0,1){2}}
\put(7.5,22.5){\line(-1,0){2}}
\put(7.5,22.5){\line(0,-1){2}}
\put(5.5,22.5){\line(1,-1){2}}
\put(6.5,20.5){\line(0,1){2}}
\put(5.5,21.5){\line(1,0){2}}
\put(7.1,21.6){$\scriptstyle B$}
\put(6.6,22.05){$\scriptstyle B$}
\put(6.1,22.05){$\scriptstyle B$}
\put(5.6,21.6){$\scriptstyle \overline{Y}$}
\put(5.6,21.05){$\scriptstyle \overline{Y}$}
\put(6.1,20.6){$\scriptstyle \overline{B}$}
\put(6.6,20.6){$\scriptstyle \overline{B}$}
\put(7.1,21.05){$\scriptstyle \overline{B}$}

\put(0,18.8){$\overline{c}$)}
\put(1,18.5){\line(1,0){1}}
\put(1,19.5){\line(1,-1){1}}
\put(1,18.5){\line(0,1){1}}
\put(2,19.5){\line(-1,0){1}}
\put(2,19.5){\line(0,-1){1}}
\put(1.1,18.6){$\scriptstyle \overline{C}$}
\put(1.6,19.05){$\scriptstyle \overline{C}$}
\put(2.15,18.55){,}
\put(2.5,18.5){\line(1,0){1}}
\put(2.5,19.5){\line(1,-1){1}}
\put(2.5,18.5){\line(0,1){1}}
\put(3.5,19.5){\line(-1,0){1}}
\put(3.5,19.5){\line(0,-1){1}}
\put(2.6,18.6){$\scriptstyle \overline{X}$}
\put(3.1,19.05){$\scriptstyle \overline{X}$}
{
\thicklines
\put(4,19){\vector(1,0){1}}
}
\put(5.5,18){\line(1,0){2}}
\put(5.5,18){\line(1,1){2}}
\put(5.5,18){\line(0,1){2}}
\put(7.5,20){\line(-1,0){2}}
\put(7.5,20){\line(0,-1){2}}
\put(5.5,20){\line(1,-1){2}}
\put(6.5,18){\line(0,1){2}}
\put(5.5,19){\line(1,0){2}}
\put(7.1,19.1){$\scriptstyle \overline{X}$}
\put(6.6,19.55){$\scriptstyle \overline{C}$}
\put(6.1,19.55){$\scriptstyle \overline{C}$}
\put(5.6,19.1){$\scriptstyle \overline{C}$}
\put(5.6,18.55){$\scriptstyle C$}
\put(6.1,18.1){$\scriptstyle C$}
\put(6.6,18.1){$\scriptstyle C$}
\put(7.1,18.55){$\scriptstyle \overline{X}$}

\put(0,16.5){$\overline{d}$)}
\put(1,16.0){\line(1,0){1}}
\put(1,17.0){\line(1,-1){1}}
\put(1,16.0){\line(0,1){1}}
\put(2,17.0){\line(-1,0){1}}
\put(2,17.0){\line(0,-1){1}}
\put(1.1,16.1){$\scriptstyle \overline{D}$}
\put(1.6,16.55){$\scriptstyle \overline{D}$}
\put(2.15,16.05){,}
\put(2.5,16.0){\line(1,0){1}}
\put(2.5,17.0){\line(1,-1){1}}
\put(2.5,16.0){\line(0,1){1}}
\put(3.5,17.0){\line(-1,0){1}}
\put(3.5,17.0){\line(0,-1){1}}
\put(2.6,16.1){$\scriptstyle \overline{V}$}
\put(3.1,16.55){$\scriptstyle \overline{V}$}
{
\thicklines
\put(4,16.5){\vector(1,0){1}}
}
\put(5.5,15.5){\line(1,0){2}}
\put(5.5,15.5){\line(1,1){2}}
\put(5.5,15.5){\line(0,1){2}}
\put(7.5,17.5){\line(-1,0){2}}
\put(7.5,17.5){\line(0,-1){2}}
\put(5.5,17.5){\line(1,-1){2}}
\put(6.5,15.5){\line(0,1){2}}
\put(5.5,16.5){\line(1,0){2}}
\put(7.1,16.6){$\scriptstyle D$}
\put(6.6,17.05){$\scriptstyle D$}
\put(6.1,17.05){$\scriptstyle \overline{D}$}
\put(5.6,16.6){$\scriptstyle \overline{D}$}
\put(5.6,16.05){$\scriptstyle \overline{D}$}
\put(6.1,15.6){$\scriptstyle \overline{V}$}
\put(6.6,15.6){$\scriptstyle \overline{V}$}
\put(7.1,16.05){$\scriptstyle D$}

\put(0,13.8){$\overline{e}$)}
\put(1,13.5){\line(1,0){1}}
\put(1,14.5){\line(1,-1){1}}
\put(1,13.5){\line(0,1){1}}
\put(2,14.5){\line(-1,0){1}}
\put(2,14.5){\line(0,-1){1}}
\put(1.1,13.6){$\scriptstyle A$}
\put(1.6,14.05){$\scriptstyle \overline{Z}$}
\put(2.15,13.55){,}
\put(2.5,13.5){\line(1,0){1}}
\put(2.5,14.5){\line(1,-1){1}}
\put(2.5,13.5){\line(0,1){1}}
\put(3.5,14.5){\line(-1,0){1}}
\put(3.5,14.5){\line(0,-1){1}}
\put(2.6,13.6){$\scriptstyle Z$}
\put(3.1,14.05){$\scriptstyle \overline{Z}$}
\put(3.65,13.55){,}
\put(4,13.5){\line(1,0){1}}
\put(4,14.5){\line(1,-1){1}}
\put(4,13.5){\line(0,1){1}}
\put(5,14.5){\line(-1,0){1}}
\put(5,14.5){\line(0,-1){1}}
\put(4.1,13.6){$\scriptstyle Z$}
\put(4.6,14.05){$\scriptstyle \overline{X}$}
{
\thicklines
\put(5.5,14){\vector(1,0){1}}
}
\put(7.0,13){\line(1,0){2}}
\put(7.0,13){\line(1,1){2}}
\put(7.0,13){\line(0,1){2}}
\put(9.0,15){\line(-1,0){2}}
\put(9.0,15){\line(0,-1){2}}
\put(7.0,15){\line(1,-1){2}}
\put(8.0,13){\line(0,1){2}}
\put(7.0,14){\line(1,0){2}}
\put(8.6,14.1){$\scriptstyle \overline{X}$}
\put(8.1,14.55){$\scriptstyle \overline{Z}$}
\put(7.6,14.55){$\scriptstyle \overline{Z}$}
\put(7.1,14.1){$\scriptstyle Z$}
\put(7.1,13.55){$\scriptstyle Z$}
\put(7.6,13.1){$\scriptstyle Z$}
\put(8.1,13.1){$\scriptstyle Z$}
\put(8.6,13.55){$\scriptstyle \overline{X}$}

\put(0,11.3){$\overline{f}$)}
\put(1,11){\line(1,0){1}}
\put(1,12){\line(1,-1){1}}
\put(1,11){\line(0,1){1}}
\put(2,12){\line(-1,0){1}}
\put(2,12){\line(0,-1){1}}
\put(1.1,11.1){$\scriptstyle \overline{Y}$}
\put(1.6,11.55){$\scriptstyle B$}
\put(2.15,11.05){,}
\put(2.5,11){\line(1,0){1}}
\put(2.5,12){\line(1,-1){1}}
\put(2.5,11){\line(0,1){1}}
\put(3.5,12){\line(-1,0){1}}
\put(3.5,12){\line(0,-1){1}}
\put(2.6,11.1){$\scriptstyle \overline{Y}$}
\put(3.1,11.55){$\scriptstyle Y$}
\put(3.65,11.05){,}
\put(4,11){\line(1,0){1}}
\put(4,12){\line(1,-1){1}}
\put(4,11){\line(0,1){1}}
\put(5,12){\line(-1,0){1}}
\put(5,12){\line(0,-1){1}}
\put(4.1,11.1){$\scriptstyle \overline{V}$}
\put(4.6,11.55){$\scriptstyle Y$}
{
\thicklines
\put(5.5,11.5){\vector(1,0){1}}
}
\put(7.0,10.5){\line(1,0){2}}
\put(7.0,10.5){\line(1,1){2}}
\put(7.0,10.5){\line(0,1){2}}
\put(9.0,12.5){\line(-1,0){2}}
\put(9.0,12.5){\line(0,-1){2}}
\put(7.0,12.5){\line(1,-1){2}}
\put(8.0,10.5){\line(0,1){2}}
\put(7.0,11.5){\line(1,0){2}}
\put(8.6,11.6){$\scriptstyle Y$}
\put(8.1,12.05){$\scriptstyle Y$}
\put(7.6,12.05){$\scriptstyle Y$}
\put(7.1,11.6){$\scriptstyle \overline{Y}$}
\put(7.1,11.05){$\scriptstyle \overline{Y}$}
\put(7.6,10.6){$\scriptstyle \overline{V}$}
\put(8.1,10.6){$\scriptstyle \overline{V}$}
\put(8.6,11.05){$\scriptstyle Y$}

\put(0,8.8){$\overline{g}$)}
\put(1,8.5){\line(1,0){1}}
\put(1,9.5){\line(1,-1){1}}
\put(1,8.5){\line(0,1){1}}
\put(2,9.5){\line(-1,0){1}}
\put(2,9.5){\line(0,-1){1}}
\put(1.1,8.6){$\scriptstyle C$}
\put(1.6,9.05){$\scriptstyle \overline{X}$}
\put(2.15,8.55){,}
\put(2.5,8.5){\line(1,0){1}}
\put(2.5,9.5){\line(1,-1){1}}
\put(2.5,8.5){\line(0,1){1}}
\put(3.5,9.5){\line(-1,0){1}}
\put(3.5,9.5){\line(0,-1){1}}
\put(2.6,8.6){$\scriptstyle X$}
\put(3.1,9.05){$\scriptstyle \overline{X}$}
\put(3.65,8.55){,}
\put(4,8.5){\line(1,0){1}}
\put(4,9.5){\line(1,-1){1}}
\put(4,8.5){\line(0,1){1}}
\put(5,9.5){\line(-1,0){1}}
\put(5,9.5){\line(0,-1){1}}
\put(4.1,8.6){$\scriptstyle X$}
\put(4.6,9.05){$\scriptstyle \overline{Z}$}
{
\thicklines
\put(5.5,9){\vector(1,0){1}}
}
\put(7.0,8){\line(1,0){2}}
\put(7.0,8){\line(1,1){2}}
\put(7.0,8){\line(0,1){2}}
\put(9.0,10){\line(-1,0){2}}
\put(9.0,10){\line(0,-1){2}}
\put(7.0,10){\line(1,-1){2}}
\put(8.0,8){\line(0,1){2}}
\put(7.0,9){\line(1,0){2}}
\put(8.6,9.1){$\scriptstyle \overline{X}$}
\put(8.1,9.55){$\scriptstyle \overline{Z}$}
\put(7.6,9.55){$\scriptstyle \overline{Z}$}
\put(7.1,9.1){$\scriptstyle X$}
\put(7.1,8.55){$\scriptstyle X$}
\put(7.6,8.1){$\scriptstyle X$}
\put(8.1,8.1){$\scriptstyle X$}
\put(8.6,8.55){$\scriptstyle \overline{X}$}

\put(0,6.3){$\overline{h}$)}
\put(1,6){\line(1,0){1}}
\put(1,7){\line(1,-1){1}}
\put(1,6){\line(0,1){1}}
\put(2,7){\line(-1,0){1}}
\put(2,7){\line(0,-1){1}}
\put(1.1,6.1){$\scriptstyle \overline{V}$}
\put(1.6,6.55){$\scriptstyle D$}
\put(2.15,6.05){,}
\put(2.5,6){\line(1,0){1}}
\put(2.5,7){\line(1,-1){1}}
\put(2.5,6){\line(0,1){1}}
\put(3.5,7){\line(-1,0){1}}
\put(3.5,7){\line(0,-1){1}}
\put(2.6,6.1){$\scriptstyle \overline{V}$}
\put(3.1,6.55){$\scriptstyle V$}
\put(3.65,6.05){,}
\put(4,6){\line(1,0){1}}
\put(4,7){\line(1,-1){1}}
\put(4,6){\line(0,1){1}}
\put(5,7){\line(-1,0){1}}
\put(5,7){\line(0,-1){1}}
\put(4.1,6.1){$\scriptstyle \overline{Y}$}
\put(4.6,6.55){$\scriptstyle V$}
{
\thicklines
\put(5.5,6.5){\vector(1,0){1}}
}
\put(7.0,5.5){\line(1,0){2}}
\put(7.0,5.5){\line(1,1){2}}
\put(7.0,5.5){\line(0,1){2}}
\put(9.0,7.5){\line(-1,0){2}}
\put(9.0,7.5){\line(0,-1){2}}
\put(7.0,7.5){\line(1,-1){2}}
\put(8.0,5.5){\line(0,1){2}}
\put(7.0,6.5){\line(1,0){2}}
\put(8.6,6.6){$\scriptstyle V$}
\put(8.1,7.05){$\scriptstyle V$}
\put(7.6,7.05){$\scriptstyle V$}
\put(7.1,6.6){$\scriptstyle \overline{Y}$}
\put(7.1,6.05){$\scriptstyle \overline{Y}$}
\put(7.6,5.6){$\scriptstyle \overline{V}$}
\put(8.1,5.6){$\scriptstyle \overline{V}$}
\put(8.6,6.05){$\scriptstyle V$}

\put(0,3.8){$\overline{i}$)}
\put(1,3.5){\line(1,0){1}}
\put(1,4.5){\line(1,-1){1}}
\put(1,3.5){\line(0,1){1}}
\put(2,4.5){\line(-1,0){1}}
\put(2,4.5){\line(0,-1){1}}
\put(1.1,3.6){$\scriptstyle Z$}
\put(1.6,4.05){$\scriptstyle Y$}
\put(2.15,3.55){,}
\put(2.5,3.5){\line(1,0){1}}
\put(2.5,4.5){\line(1,-1){1}}
\put(2.5,3.5){\line(0,1){1}}
\put(3.5,4.5){\line(-1,0){1}}
\put(3.5,4.5){\line(0,-1){1}}
\put(2.6,3.6){$\scriptstyle Y$}
\put(3.1,4.05){$\scriptstyle Z$}
\put(3.65,3.55){,}
\put(4,3.5){\line(1,0){1}}
\put(4,4.5){\line(1,-1){1}}
\put(4,3.5){\line(0,1){1}}
\put(5,4.5){\line(-1,0){1}}
\put(5,4.5){\line(0,-1){1}}
\put(4.1,3.6){$\scriptstyle A$}
\put(4.6,4.05){$\scriptstyle Z$}
\put(5.15,3.55){,}
\put(5.5,3.5){\line(1,0){1}}
\put(5.5,4.5){\line(1,-1){1}}
\put(5.5,3.5){\line(0,1){1}}
\put(6.5,4.5){\line(-1,0){1}}
\put(6.5,4.5){\line(0,-1){1}}
\put(5.6,3.6){$\scriptstyle Y$}
\put(6.1,4.05){$\scriptstyle B$}
{
\thicklines
\put(7.0,4){\vector(1,0){1}}
}
\put(8.5,3){\line(1,0){2}}
\put(8.5,3){\line(1,1){2}}
\put(8.5,3){\line(0,1){2}}
\put(10.5,5){\line(-1,0){2}}
\put(10.5,5){\line(0,-1){2}}
\put(8.5,5){\line(1,-1){2}}
\put(9.5,3){\line(0,1){2}}
\put(8.5,4){\line(1,0){2}}
\put(10.1,4.1){$\scriptstyle Y$}
\put(9.6,4.55){$\scriptstyle Y$}
\put(9.1,4.55){$\scriptstyle Y$}
\put(8.6,4.1){$\scriptstyle Z$}
\put(8.6,3.55){$\scriptstyle Z$}
\put(9.1,3.1){$\scriptstyle Z$}
\put(9.6,3.1){$\scriptstyle Z$}
\put(10.1,3.55){$\scriptstyle Y$}

\put(0,1.3){$\overline{k}$)}
\put(1,1){\line(1,0){1}}
\put(1,2){\line(1,-1){1}}
\put(1,1){\line(0,1){1}}
\put(2,2){\line(-1,0){1}}
\put(2,2){\line(0,-1){1}}
\put(1.1,1.1){$\scriptstyle X$}
\put(1.6,1.55){$\scriptstyle V$}
\put(2.15,1.05){,}
\put(2.5,1){\line(1,0){1}}
\put(2.5,2){\line(1,-1){1}}
\put(2.5,1){\line(0,1){1}}
\put(3.5,2){\line(-1,0){1}}
\put(3.5,2){\line(0,-1){1}}
\put(2.6,1.1){$\scriptstyle V$}
\put(3.1,1.55){$\scriptstyle X$}
\put(3.65,1.05){,}
\put(4,1){\line(1,0){1}}
\put(4,2){\line(1,-1){1}}
\put(4,1){\line(0,1){1}}
\put(5,2){\line(-1,0){1}}
\put(5,2){\line(0,-1){1}}
\put(4.1,1.1){$\scriptstyle V$}
\put(4.6,1.55){$\scriptstyle D$}
\put(5.15,1.05){,}
\put(5.5,1){\line(1,0){1}}
\put(5.5,2){\line(1,-1){1}}
\put(5.5,1){\line(0,1){1}}
\put(6.5,2){\line(-1,0){1}}
\put(6.5,2){\line(0,-1){1}}
\put(5.6,1.1){$\scriptstyle C$}
\put(6.1,1.55){$\scriptstyle X$}
{
\thicklines
\put(7.0,1.5){\vector(1,0){1}}
}
\put(8.5,0.5){\line(1,0){2}}
\put(8.5,0.5){\line(1,1){2}}
\put(8.5,0.5){\line(0,1){2}}
\put(10.5,2.5){\line(-1,0){2}}
\put(10.5,2.5){\line(0,-1){2}}
\put(8.5,2.5){\line(1,-1){2}}
\put(9.5,0.5){\line(0,1){2}}
\put(8.5,1.5){\line(1,0){2}}
\put(10.1,1.6){$\scriptstyle V$}
\put(9.6,2.05){$\scriptstyle V$}
\put(9.1,2.05){$\scriptstyle V$}
\put(8.6,1.6){$\scriptstyle X$}
\put(8.6,1.05){$\scriptstyle X$}
\put(9.1,0.6){$\scriptstyle X$}
\put(9.6,0.6){$\scriptstyle X$}
\put(10.1,1.05){$\scriptstyle V$}
\end{picture}
\end{tabular}
\caption{Substitution on squares, derived from the substitution in
  Fig.~\ref{fig:haussub}.   \label{fig:sqsub}}
\end{figure}
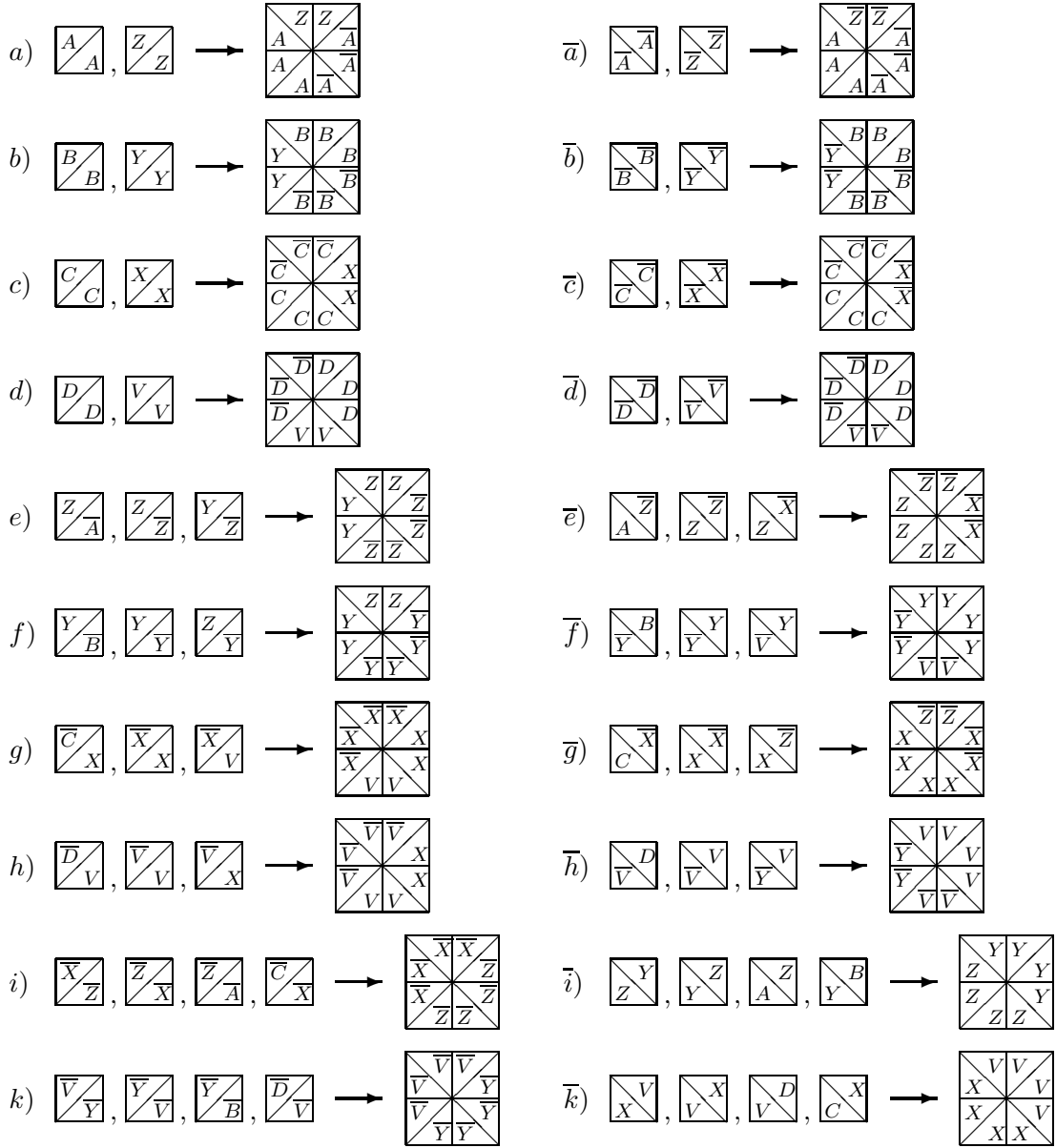

We then obtain the substitution for ``squares'' in
Fig.~\ref{fig:sqsub}, p.~\pageref{fig:sqsub}, which is again a
substitution on $L=\Z^2_{}$ by a factor of $2$.

\begin{figure} 
\begin{center}
\setlength{\unitlength}{.7cm}
\begin{picture}(22,21)
\put(3.5,18){\oval(4,1)}
\put(3.5,18){\oval(4.2,1.2)}
\put(1.7,17.9){$\scriptstyle \{a,b,c,d,e,f,g,h,i,k\}$}

\put(8.5,18){\oval(3.4,1)}
\put(7,17.9){$\scriptstyle \{b,d,e,f,g,h,i,k\}$}

\put(13.5,18){\oval(2.6,1)}
\put(12.4,17.9){$\scriptstyle \{e,f,g,h,i,k\}$}

\put(18.5,18){\oval(3.4,1)}
\put(17,17.9){$\scriptstyle \{a,c,e,f,g,h,i,k\}$}

\put(8,15){\oval(1.7,1)}
\put(7.35,14.9){$\scriptstyle \{\overline{a},\overline{b},\overline{k}\}$}

\put(14,15){\oval(1.7,1)}
\put(13.4,14.9){$\scriptstyle \{\overline{c},\overline{d},\overline{i}\}$}

\put(2,12){\oval(1.7,1)}
\put(1.35,11.9){$\scriptstyle \{\overline{e},\overline{f},\overline{k}\}$}

\put(20,12){\oval(1.7,1)}
\put(19.4,11.9){$\scriptstyle \{\overline{g},\overline{h},\overline{i}\}$}

\put(2,9){\oval(1.7,1)}
\put(1.35,8.9){$\scriptstyle \{e,g,k\}$}

\put(20,9){\oval(1.7,1)}
\put(19.4,8.9){$\scriptstyle \{f,h,i\}$}

\put(8,6){\oval(1.7,1)}
\put(7.35,5.9){$\scriptstyle \{a,c,k\}$}

\put(14,6){\oval(1.7,1)}
\put(13.4,5.9){$\scriptstyle \{b,d,i\}$}

\put(3.5,3){\oval(4,1)}
\put(3.5,3){\oval(4.2,1.2)}
\put(1.7,2.9){$\scriptstyle
  \{\overline{a},\overline{b},\overline{c},\overline{d},\overline{e},\overline{f},\overline{g},\overline{h},\overline{i},\overline{k}\}$} 

\put(8.5,3){\oval(3.4,1)}
\put(7,2.9){$\scriptstyle
  \{\overline{c},\overline{d},\overline{e},\overline{f},\overline{g},\overline{h},\overline{i},\overline{k}\}$} 

\put(13.5,3){\oval(2.6,1)}
\put(12.4,2.9){$\scriptstyle
  \{\overline{e},\overline{f},\overline{g},\overline{h},\overline{i},\overline{k}\}$} 

\put(18.5,3){\oval(3.4,1)}
\put(17,2.9){$\scriptstyle
  \{\overline{a},\overline{b},\overline{e},\overline{f},\overline{g},\overline{h},\overline{i},\overline{k}\}$} 
\qbezier[1000](3.5,18.7)(11,22)(18.5,18.6)
\put(11.1,20.32){\vector(1,0){0}}

\qbezier(8.2,18.6)(8,18.8)(8,19)
\qbezier(8.5,19.5)(8,19.5)(8,19)
\qbezier(8.5,19.5)(9,19.5)(9,19)
\qbezier(8.8,18.6)(9,18.8)(9,19)
\put(8.65,19.5){\vector(1,0){0}}

\qbezier(13.2,18.6)(13,18.8)(13,19)
\qbezier(13.5,19.5)(13,19.5)(13,19)
\qbezier(13.5,19.5)(14,19.5)(14,19)
\qbezier(13.8,18.6)(14,18.8)(14,19)
\put(13.65,19.5){\vector(1,0){0}}

\put(5.8,18){\vector(1,0){.8}}
\put(10.4,18){\vector(1,0){1.6}}
\put(16.6,18){\vector(-1,0){1.6}}

\qbezier(20.3,18.3)(20.5,18.5)(20.7,18.5)
\qbezier(21.2,18)(21.2,18.5)(20.7,18.5)
\qbezier(21.2,18)(21.2,17.5)(20.7,17.5)
\qbezier(20.3,17.7)(20.5,17.5)(20.7,17.5)
\put(21.2,17.9){\vector(0,-1){0}}

\put(4.1,17.2){\vector(2,-1){3.2}}
\put(8,17.4){\vector(0,-1){1.8}}
\put(12.3,17.4){\vector(-2,-1){3.6}}
\put(17.1,17.4){\vector(-4,-1){8.2}}
\put(5.0,17.2){\vector(4,-1){8.0}}
\put(9.7,17.4){\vector(2,-1){3.6}}
\put(14,17.4){\vector(0,-1){1.8}}
\put(18.3,17.4){\vector(-2,-1){3.6}}

\qbezier(7,15.3)(6.8,15.5)(6.6,15.5)
\qbezier(6.1,15)(6.1,15.5)(6.6,15.5)
\qbezier(6.1,15)(6.1,14.5)(6.6,14.5)
\qbezier(7,14.7)(6.8,14.5)(6.6,14.5)
\put(6.1,14.9){\vector(0,-1){0}}

\qbezier(15,15.3)(15.2,15.5)(15.4,15.5)
\qbezier(15.9,15)(15.9,15.5)(15.4,15.5)
\qbezier(15.9,15)(15.9,14.5)(15.4,14.5)
\qbezier(15,14.7)(15.2,14.5)(15.4,14.5)
\put(15.9,14.9){\vector(0,-1){0}}

\put(7.1,14.4){\vector(-2,-1){4.1}}
\put(14.9,14.4){\vector(2,-1){4.1}}
\put(2.6,11.4){\vector(1,-1){4.8}}
\put(2.8,11.6){\vector(2,-1){10.3}}
\put(19.4,11.4){\vector(-1,-1){4.8}}
\put(19.2,11.6){\vector(-2,-1){10.3}}

\qbezier(1,12.3)(.8,12.5)(.6,12.5)
\qbezier(.1,12)(.1,12.5)(.6,12.5)
\qbezier(.1,12)(.1,11.5)(.6,11.5)
\qbezier(1,11.7)(.8,11.5)(.6,11.5)
\put(.1,11.9){\vector(0,-1){0}}

\qbezier(21,12.3)(21.2,12.5)(21.4,12.5)
\qbezier(21.9,12)(21.9,12.5)(21.4,12.5)
\qbezier(21.9,12)(21.9,11.5)(21.4,11.5)
\qbezier(21,11.7)(21.2,11.5)(21.4,11.5)
\put(21.9,11.9){\vector(0,-1){0}}

\put(7.8,14.4){\vector(0,-1){7.8}}
\put(8,6.6){\vector(0,1){7.8}}

\put(8.3,14.4){\vector(2,-3){5.2}}
\put(13.7,6.6){\vector(-2,3){5.2}}

\put(13.5,14.4){\vector(-2,-3){5.2}}
\put(8.5,6.6){\vector(2,3){5.2}}

\put(14,14.4){\vector(0,-1){7.8}}
\put(14.2,6.6){\vector(0,1){7.8}}

\qbezier(1,9.3)(.8,9.5)(.6,9.5)
\qbezier(.1,9)(.1,9.5)(.6,9.5)
\qbezier(.1,9)(.1,8.5)(.6,8.5)
\qbezier(1,8.7)(.8,8.5)(.6,8.5)
\put(.1,8.9){\vector(0,-1){0}}

\qbezier(21,9.3)(21.2,9.5)(21.4,9.5)
\qbezier(21.9,9)(21.9,9.5)(21.4,9.5)
\qbezier(21.9,9)(21.9,8.5)(21.4,8.5)
\qbezier(21,8.7)(21.2,8.5)(21.4,8.5)
\put(21.9,8.9){\vector(0,-1){0}}

\put(7.1,6.6){\vector(-2,1){4.1}}
\put(14.9,6.6){\vector(2,1){4.1}}
\put(2.6,9.6){\vector(1,1){4.8}}
\put(2.8,9.4){\vector(2,1){10.3}}
\put(19.4,9.6){\vector(-1,1){4.8}}
\put(19.2,9.4){\vector(-2,1){10.3}}

\qbezier(7,6.3)(6.8,6.5)(6.6,6.5)
\qbezier(6.1,6)(6.1,6.5)(6.6,6.5)
\qbezier(6.1,6)(6.1,5.5)(6.6,5.5)
\qbezier(7,5.7)(6.8,5.5)(6.6,5.5)
\put(6.1,5.9){\vector(0,-1){0}}

\qbezier(15,6.3)(15.2,6.5)(15.4,6.5)
\qbezier(15.9,6)(15.9,6.5)(15.4,6.5)
\qbezier(15.9,6)(15.9,5.5)(15.4,5.5)
\qbezier(15,5.7)(15.2,5.5)(15.4,5.5)
\put(15.9,5.9){\vector(0,-1){0}}

\put(4.1,3.8){\vector(2,1){3.2}}
\put(8,3.6){\vector(0,1){1.8}}
\put(12.3,3.6){\vector(-2,1){3.6}}
\put(17.1,3.6){\vector(-4,1){8.2}}
\put(5.0,3.8){\vector(4,1){8.0}}
\put(9.7,3.6){\vector(2,1){3.6}}
\put(14,3.6){\vector(0,1){1.8}}
\put(18.3,3.6){\vector(-2,1){3.6}}

\put(5.8,3){\vector(1,0){.8}}
\put(10.4,3){\vector(1,0){1.6}}
\put(16.6,3){\vector(-1,0){1.6}}

\qbezier(20.3,3.3)(20.5,3.5)(20.7,3.5)
\qbezier(21.2,3)(21.2,3.5)(20.7,3.5)
\qbezier(21.2,3)(21.2,2.5)(20.7,2.5)
\qbezier(20.3,2.7)(20.5,2.5)(20.7,2.5)
\put(21.2,2.9){\vector(0,-1){0}}

\qbezier(8.2,2.4)(8,2.2)(8,2)
\qbezier(8.5,1.5)(8,1.5)(8,2)
\qbezier(8.5,1.5)(9,1.5)(9,2)
\qbezier(8.8,2.4)(9,2.2)(9,2)
\put(8.65,1.5){\vector(1,0){0}}

\qbezier(13.2,2.4)(13,2.2)(13,2)
\qbezier(13.5,1.5)(13,1.5)(13,2)
\qbezier(13.5,1.5)(14,1.5)(14,2)
\qbezier(13.8,2.4)(14,2.2)(14,2)
\put(13.65,1.5){\vector(1,0){0}}

\qbezier[1000](3.5,2.3)(11,-1)(18.5,2.4)
\put(11.1,.67){\vector(1,0){0}}
\end{picture}
\end{center}
\caption{The coincidence graph of the derived tiling of the
  Semi-Detached House Tiling.
 \label{fig:hausgraph}}
\end{figure}
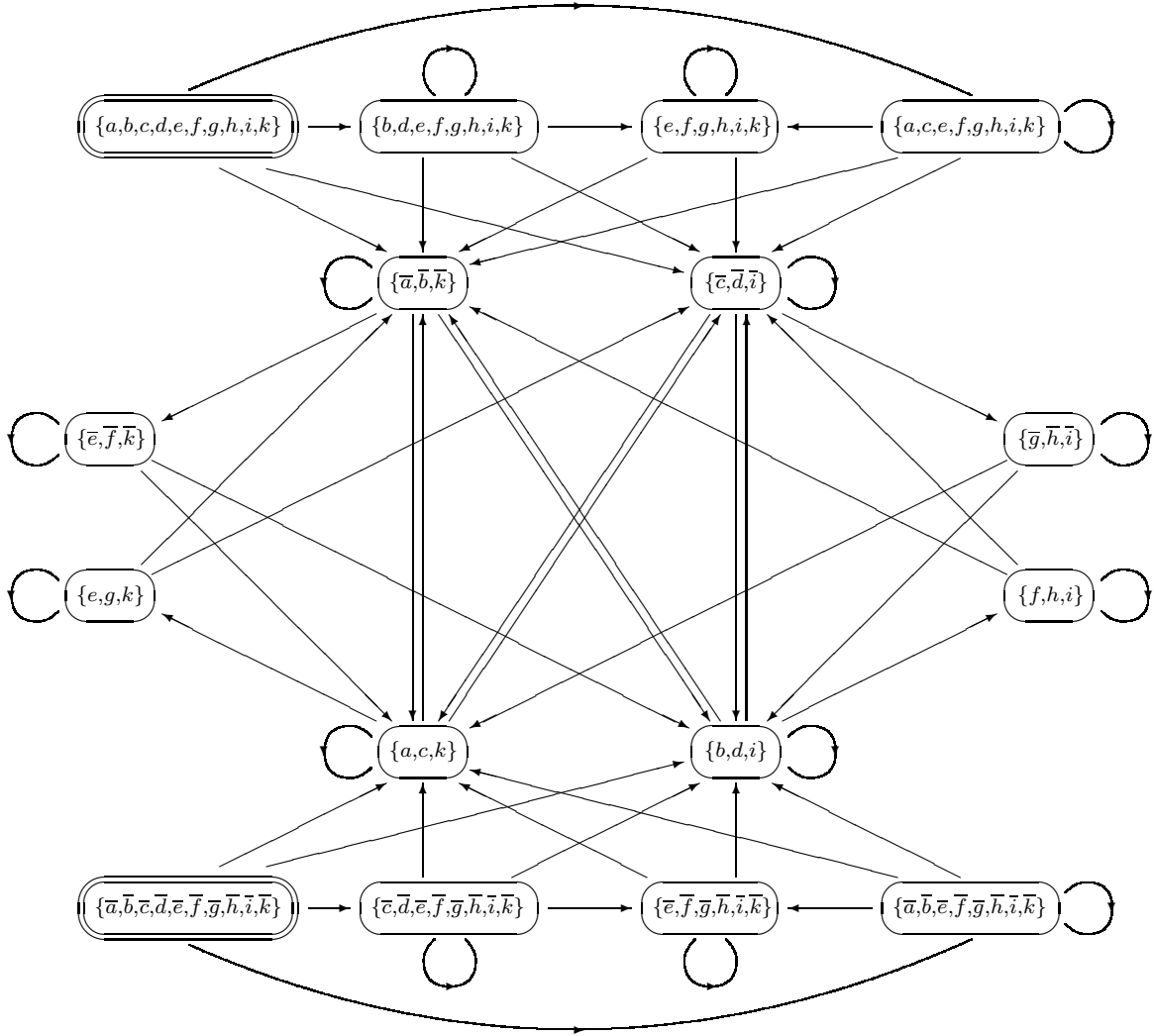

Instead of using all $56$ squares which occur in this substitution, we
form $20$ classes which are already present in Fig.~\ref{fig:sqsub} (i.e.,
squares which yield the same substitute are in the same class). We now
look at the substitution induced on these classes, e.g., 
\begin{equation*}
a \to \begin{array}{cc} \overline{i} & e \\ a & \overline{a}
\end{array} \qquad \qquad
\overline{a} \to \begin{array}{cc} \overline{e} & i \\ a & \overline{a}
\end{array}.
\end{equation*}
Obviously, if we do not have a modular coincidence for this
induced substitution, then there is no coincidence for the
original substitution either. The coincidence graph for the induced
substitution is shown in Fig.~\ref{fig:hausgraph},
p.~\pageref{fig:hausgraph}. As suggested by Fig.~\ref{fig:sqsub},
$L'$ forms a checkerboard with respect to $L$, i.e.,
$[L:L']=2$. 

As we do not have a modular coincidence for this substitution (no vertex
in the coincidence graph represents a singleton), the
semi-detached house tiling does not arise from a model set and is not
pure point diffractive. This non-trivial example shows that
analyzing LSS and tilings can often be quite elaborate.   

\begin{appendix}

\section{Modular Coincidence vs. Dekking Coincidence}

An admissible LSS in one dimension (i.e., on
$L=\Z$) is known as a primitive \emph{substitution of constant length}
$q$, see~\cite{Dek78,Queff,Fogg}. Here the mappings of the MFS $\Phi$
are of the form 
$x\mapsto q\cdot x+a^{}_{ijk}$ where $a^{}_{ijk}\in\{0,\ldots,q-1\}$,
i.e., each letter is substituted by exactly $q$ letters. 

The first step in analyzing such a substitution is to determine
$L'$. Since $L=\Z$, we have $L'=r\cdot\Z$ and $r\in\N$ is to be
determined. To this end, we define, for a sequence $u=\ldots
u^{}_{-2}u^{}_{-1}u^{}_{0}u^{}_{1}u^{}_{2}\ldots$, the sets
$S^{}_k=\{j\in\Z\mathbin| u^{}_{j+k}=u^{}_k\}$ and numbers
$g^{}_k=\gcd S^{}_k$. It is easy to see that $g^{}_k=g^{}_{k'}$ if
$u^{}_k=u^{}_{k'}$, so that we only have to determine the $g^{}_k$
once for each letter of the alphabet. Furthermore, we have $L^{}_i
= g^{}_k\cdot\Z$ if $u^{}_k=i$ (in other words: if $u_k$ is of type
$i$) and, by $L'=L^{}_1+\ldots +L^{}_m$ ($m$
is the cardinality of the alphabet), also $r=\gcd\{g^{}_k\mathbin| k\in\Z\}$. 

Let us introduce the \emph{height} $h$ of a substitution of constant
length $q$ (see~\cite[Definition II.8]{Dek78} and~\cite[Definition VI.1.]{Queff}):
$h=\max\{n\ge 1\mathbin| \gcd\{n,q\}=1 \textnormal{ and } n \textnormal{
  divides } g^{}_0\}$. We have the following property for every
$k\in\Z$ (\cite[Remark II.9(i)]{Dek78} and~\cite[VI.2.1.1.]{Queff}): 
\begin{equation*}
\{n\ge 1\mathbin| \gcd\{n,q\}=1 \textnormal{ and } n \textnormal{
  divides } g^{}_0\} = \{n\ge 1\mathbin| \gcd\{n,q\}=1 \textnormal{
  and } n \textnormal{ divides } g^{}_k\}
\end{equation*}    
From this, it follows directly that $1\le h\le r\le m$ and $h$ divides
$r$, but is relatively prime to $q$. If $r=m$, then the sequence under
consideration is periodic (compare to the argument in the proof of
Corollary~\ref{cor:bij}). It can happen that $h<r$. 

\smallskip

\textsc{Example:} The paperfolding sequence can be generated by the
substitution $a\mapsto ab$, $b\mapsto cb$, $c\mapsto ad$, $d\mapsto
cd$ and a subsequent projection $\{a,b\}\mapsto 0$, $\{c,d\}\mapsto
1$ on a two letter alphabet $\{0,1\}$. We consider the original
substitution on $\{a,b,c,d\}$ for which we have: $q=2$, 
$g^{}_k=2$ for all $k\in\Z$, $r=2$ and $h=1$. 

\smallskip

For a substitution $\sigma$ of constant length with $h>1$, it is always
possible to introduce a substitution $\eta$ on a (usually bigger) alphabet of
the same constant length with height $1$ such that the sequence
generated by $\sigma$ can also be obtained from the sequence generated
by $\eta$, see~\cite[Lemma II.17 \& Definition II.20]{Dek78}. We illustrate the
general method by an example, also compare this construction to
Theorem~\ref{thm:MLDtiling}.   
 
\smallskip

\textsc{Example:} Consider the substitution (see~\cite[Example
  following Definition II.20]{Dek78}) \linebreak $0\mapsto 010$, $1\mapsto
  102$ and $2\mapsto 201$ of height $2$. This substitution generates a sequence
  $u=\ldots 010102010102\ldots$. Now we take as the new alphabet all
  words in $\{u^{}_{n\cdot h}\ldots u^{}_{n\cdot h+h-1}\mathbin|
  n\in\Z\}$, in this case just $a=01$ and $b=02$. Then the
  substitution $\eta$ of height $1$ is given by $a\mapsto aab$ and
  $b\mapsto aba$ (since $a=01\mapsto 010102=aab$). 

\smallskip

Let $\eta$ be a substitution of constant length $q$ and height
$h=1$. We say $\eta$ admits a \emph{coincidence} (or
\emph{Dekking coincidence}) if there exist $k$ and $0\le j< q_{}^k$ such
that the $j$-th symbol of $\eta^k_{}(i)$ is the same for all letters
$i$ (\cite[Definition VI.3]{Queff}, see also~\cite[Defnition III.1]{Dek78}). 
In the above notation, $\eta$ admits a coincidence if
$\eta^{k}_{}(i)^{}_j$ is the same for all $i$. A
substitution $\sigma$ of constant length $q$ and height $h>1$ is said
to admit a coincidence if its associated substitution of height $h=1$
admits one. 

Reformulating~\cite[Theorem II.13 \& Theorem III.7]{Dek78} (compare
to~\cite[Theorem VI.13 \& Theorem VI.24]{Queff} and~\cite[Theorem
7.3.1 \& Theorem 7.3.6]{Fogg}) we obtain:

\begin{lemma}
Let $\sigma$ be a substitution of constant length $q$ and height
$h$. Let $u$ be the sequence generated by $\sigma$. Then $u$ is a
regular model set iff $\sigma$ admits a Dekking coincidence. The
internal space $H$ is given by $H=\Z^{}_q \times C^{}_h$, where
$\Z^{}_q$ denotes the product over the distinct primes $p$ dividing $q$ of the
$p$-adic integers $\Z^{}_p$ and $C^{}_h$ is the cyclic group of order $h$.
\qed
\end{lemma}

We compare this result with the ones involving the modular coincidence
(Definition~\ref{def:modcoinc} and Theorem~\ref{thm:modcoinc}): The
substitution has a modular coincidence iff it is a regular model set
with internal space $C^{}_r\times \lim_{\gets k}
\Z/q^k_{}\Z$. Obviously, $\lim_{\gets k} \Z/q^k_{}\Z \simeq \Z^{}_q$
(e.g., observe that $\lim_{\gets k} \Z/4^k_{}\Z$ is topologically
isomorphic to $\Z^{}_2$, wherefore it is enough to consider just the distinct
primes $p$). We can write $r=h\cdot \tilde{q}$, where
$\gcd\{r,q\}=\tilde{q}$, and we observe that $C^{}_r\times\Z^{}_q$ is
topologically isomorphic to $C^{}_h\times\Z^{}_q$ (note that 
$C^{}_{\tilde{q}}\times\Z^{}_q$ and $\Z^{}_q$ are topologically isomorphic). So
both methods actually give -- up to isomorphism -- the same internal space.

In a sense, the height $h$ tells us whether $L'=r\cdot\Z$ lines up with the
length $q$ of the substitution. If $h=1$ but $r>1$, then we do not
have to know the sets $\Phi^{}_0[\cdot]$ explicitly for the
Dekking coincidence, the substitution rule itself already gives the
right partition of the alphabet. 

\smallskip

\textsc{Example continued:} The underlying substitution of the
paperfolding sequence ($a\mapsto ab$, $b\mapsto cb$, $c\mapsto ad$, $d\mapsto
cd$) has height $h=1$ but $r=2$. We observe that there is  either an
$a$ or a $c$ on the zeroth place of the substitute, while on the first
we have either $b$ or $d$. But we also have $\Psi^{}_0[0]=\{a,c\}$ and
$\Psi^{}_0[1]=\{b,d\}$, therefore the substitution lines up with this
partition. We note that this substitution has a modular coincidence
relative to $2^2\,\Z$ (after two substitutions the first letter is
always an $a$, $a\mapsto ab \mapsto abcb$, $b\mapsto cb \mapsto adcb$
etc.), and therefore is pure point diffractive.

\section{Subdivision Graph for Overlaps}

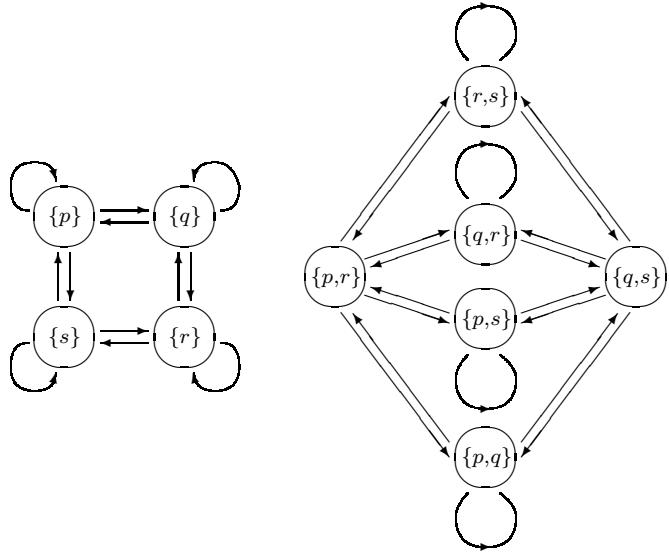
\begin{figure}[ht]
\begin{center}
\setlength{\unitlength}{.8cm}
\begin{picture}(12,10)
\put(1,6){\oval(1,1)}
\put(.73,5.9){$\scriptstyle \{p\}$}
\qbezier(.4,6.1)(.1,6.1)(.1,6.5)
\qbezier(.1,6.5)(.1,6.9)(.5,6.9)
\qbezier(.5,6.9)(.85,6.9)(.85,6.6)
\put(.86,6.6){\vector(1,-4){0}}

\put(3,6){\oval(1,1)}
\put(2.73,5.9){$\scriptstyle \{q\}$}
\qbezier(3.6,6.1)(3.9,6.1)(3.9,6.5)
\qbezier(3.9,6.5)(3.8,6.9)(3.5,6.9)
\qbezier(3.5,6.9)(3.15,6.9)(3.15,6.6)
\put(3.14,6.6){\vector(-1,-4){0}}

\put(1,4){\oval(1,1)}
\put(.73,3.9){$\scriptstyle \{s\}$}
\qbezier(.4,3.9)(.1,3.9)(.1,3.5)
\qbezier(.1,3.5)(.1,3.1)(.5,3.1)
\qbezier(.5,3.1)(.85,3.1)(.85,3.4)
\put(.86,3.4){\vector(1,4){0}}

\put(3,4){\oval(1,1)}
\put(2.73,3.9){$\scriptstyle \{r\}$}
\qbezier(3.6,3.9)(3.9,3.9)(3.9,3.5)
\qbezier(3.9,3.5)(3.9,3.1)(3.5,3.1)
\qbezier(3.5,3.1)(3.15,3.1)(3.15,3.4)
\put(3.14,3.4){\vector(-1,4){0}}

\put(1.6,6.1){\vector(1,0){.8}}
\put(2.4,5.9){\vector(-1,0){.8}}
\put(1.6,4.1){\vector(1,0){.8}}
\put(2.4,3.9){\vector(-1,0){.8}}
\put(.9,4.6){\vector(0,1){.8}}
\put(1.1,5.4){\vector(0,-1){.8}}
\put(2.9,4.6){\vector(0,1){.8}}
\put(3.1,5.4){\vector(0,-1){.8}}

\put(5.5,5){\oval(1,1)}
\put(5.1,4.9){$\scriptstyle \{p,r\}$}

\put(5.6,5.6){\vector(3,4){1.8}}
\put(7.4,7.75){\vector(-3,-4){1.6}}
\put(6,5.3){\vector(3,1){1.4}}
\put(7.4,5.6){\vector(-3,-1){1.3}}
\put(6,4.7){\vector(3,-1){1.4}}
\put(7.4,4.4){\vector(-3,1){1.3}}
\put(5.6,4.4){\vector(3,-4){1.8}}
\put(7.4,2.25){\vector(-3,4){1.6}}

\put(8,2){\oval(1,1)}
\put(7.6,1.9){$\scriptstyle \{p,q\}$}
\qbezier(7.7,1.4)(7.5,1.2)(7.5,1)
\qbezier(8.0,0.5)(7.5,0.5)(7.5,1)
\qbezier(8.0,0.5)(8.5,0.5)(8.5,1)
\qbezier(8.23,1.4)(8.5,1.2)(8.5,1)
\put(8.1,0.5){\vector(1,0){0}}

\put(8,4.3){\oval(1,1)}
\put(7.6,4.2){$\scriptstyle \{p,s\}$}
\qbezier(7.7,3.7)(7.5,3.5)(7.5,3.3)
\qbezier(8.0,2.8)(7.5,2.8)(7.5,3.3)
\qbezier(8.0,2.8)(8.5,2.8)(8.5,3.3)
\qbezier(8.23,3.7)(8.5,3.5)(8.5,3.3)
\put(8.1,2.8){\vector(1,0){0}}

\put(8,5.7){\oval(1,1)}
\put(7.6,5.6){$\scriptstyle \{q,r\}$}
\qbezier(7.7,6.3)(7.5,6.5)(7.5,6.7)
\qbezier(8.0,7.2)(7.5,7.2)(7.5,6.7)
\qbezier(8.0,7.2)(8.5,7.2)(8.5,6.7)
\qbezier(8.23,6.3)(8.5,6.5)(8.5,6.7)
\put(8.1,7.2){\vector(1,0){0}}

\put(8,8){\oval(1,1)}
\put(7.6,7.9){$\scriptstyle \{r,s\}$}
\qbezier(7.7,8.6)(7.5,8.8)(7.5,9)
\qbezier(8.0,9.5)(7.5,9.5)(7.5,9)
\qbezier(8.0,9.5)(8.5,9.5)(8.5,9)
\qbezier(8.23,8.6)(8.5,8.8)(8.5,9)
\put(8.1,9.5){\vector(1,0){0}}

\put(10.5,5){\oval(1,1)}
\put(10.1,4.9){$\scriptstyle \{q,s\}$}

\put(10.4,5.6){\vector(-3,4){1.8}}
\put(8.6,7.75){\vector(3,-4){1.6}}
\put(10,5.3){\vector(-3,1){1.4}}
\put(8.6,5.6){\vector(3,-1){1.3}}
\put(10,4.7){\vector(-3,-1){1.4}}
\put(8.6,4.4){\vector(3,1){1.3}}
\put(10.4,4.4){\vector(-3,-4){1.8}}
\put(8.6,2.25){\vector(3,4){1.6}}

\end{picture}\\
\end{center}
\caption{ \label{fig:pcoinctable}
Pair coincidence graph for the Table Tiling} 
\end{figure}

Instead of the coincidence graph $G^{}_{(\bs{V},\Phi)}$, we can also
define a \emph{pair coincidence graph} $G^{(2)}_{(\V,\Phi)}$:
$G^{(2)}_{(\V,\Phi)}$ is a directed (multi-)graph, where each set of
pairs of elements of $\Psi^{}_0[a]$ for some $a$ (i.e., both elements
are from $\Psi^{}_0[a]$ with the same $a$) is represented by a
vertex. We call pairs of the same element $\{i,i\}=\{i\}$ a
\emph{coincidence} in the pair coincidence graph. We draw an edge from
$v^{}_1$ to $v^{}_2$ iff $v^{}_2$ is a ``child''
of $v^{}_1$ (Here we mean by an ``child'': If
$v^{}_1=\{s^{}_1,\ldots,s^{}_{\ell}\}$ then
$v^{}_2=\{\Phi(s^{}_1)^{}_z,\ldots,\Phi(s^{}_{\ell})^{}_z\}$ for
some $z$.)  The pair coincidence graph of the table tiling is shown
in Fig.~\ref{fig:pcoinctable}.

Obviously, if the coincidence graph $G^{}_{(\bs{V},\Phi)}$ has a
singleton set (i.e., we have a modular coincidence), then from each vertex
of $G^{(2)}_{(\bs{V},\Phi)}$ there is a path leading to a
coincidence. Conversely, since the alphabet is finite, if from each
vertex of $G^{(2)}_{(\bs{V},\Phi)}$ there is a path leading to a
coincidence, then $G^{}_{(\bs{V},\Phi)}$ also has a singleton
set (compare this to Corollary~\ref{cor:paircoinc}). Therefore, we can
also use the pair coincidence graph to decide the 
question about modular coincidences (the pair coincidence graph of the
table tiling is not connected, the coincidences forming one connected
component; the table tiling has no modular coincidences).

\begin{figure}
\begin{center}
\setlength{\unitlength}{.8cm}
\begin{picture}(10,4)
\put(2,1.5){\oval(2,2)}
\put(1.5,1.25){\rule{.8cm}{.4cm}}
\qbezier(1.7,2.6)(1.5,2.8)(1.5,3)
\qbezier(2.0,3.5)(1.5,3.5)(1.5,3)
\qbezier(2.0,3.5)(2.5,3.5)(2.5,3)
\qbezier(2.23,2.6)(2.5,2.8)(2.5,3)
\put(2.1,3.5){\vector(1,0){0}}

\put(5.5,1.5){\oval(2,2)}
\put(5.25,1.25){\rule{.4cm}{.4cm}}
\put(4.75,1.25){\line(1,0){1.5}}
\put(4.75,1.25){\line(0,1){.5}}
\put(6.25,1.75){\line(0,-1){.5}}
\put(6.25,1.75){\line(-1,0){1.5}}

\put(6.7,1.4){\vector(1,0){1.1}}
\put(7.8,1.6){\vector(-1,0){1.1}}

\put(9,1.5){\oval(2,2)}
\put(8.5,1.5){\rule{.4cm}{.4cm}}
\put(8.5,2){\line(1,0){1}}
\put(8.5,2){\line(0,-1){1}}
\put(8.5,2){\line(0,-1){1}}
\put(9,1){\line(-1,0){.5}}
\put(9,1){\line(0,1){1}}
\put(9.5,1.5){\line(0,1){.5}}
\put(9.5,1.5){\line(-1,0){1}}

\qbezier(8.7,2.6)(8.5,2.8)(8.5,3)
\qbezier(9.0,3.5)(8.5,3.5)(8.5,3)
\qbezier(9.0,3.5)(9.5,3.5)(9.5,3)
\qbezier(9.23,2.6)(9.5,2.8)(9.5,3)
\put(9.1,3.5){\vector(1,0){0}}

\end{picture}\\
\end{center}
\caption{ \label{fig:subdivt}
Subdivision graph for overlaps for the Table Tiling} 
\end{figure}
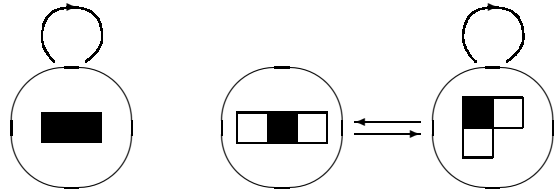

On the tiling level, the pair coincidence graph is called \emph{subdivision
graph for overlaps}, see~\cite[p.~721]{Sol97} and \cite[Appendix
2]{LMS03}. The construction of this subdivision graph follows closely,
with the necessary changes, the construction of the pair coincidence
graph. For the pair coincidence graph we can think of placing two
multi-component Delone sets on top of each other (or, placing a
translated (by a vector of $L'$) copy of the multi-component set on
top) and observing what happens under substitution. Likewise, we can
think of placing two tilings on top of each other and than we observe,
what happens to the overlaps of the tiles under substitution. As an
illustration we give the subdivision graph for overlaps (overlaps are
denoted by filled areas) for the table tiling in Fig.~\ref{fig:subdivt}
(where we identify rotated versions of the configurations), see also 
\cite[Fig.~7.4 (d)]{Sol97}. This graph should be compared with the
pair coincidence graph in Fig.~\ref{fig:pcoinctable} 
(basically, rotation identifies all the
singletons, $\{p,r\}$ and $\{q,s\}$, and the four sets $\{r,s\}$,
$\{q,r\}$, $\{p,s\}$ and $\{p,q\}$). 

\section{Connections to Automatic Sequences}

In~\cite{Haeseler}, von Haeseler is interested in the question, under which
conditions a sequence is automatic and/or generated by a (not
necessarily primitive) substitution of constant length (and in the
higher dimensional case even whether there exists an LSS).
To this end, he introduces the \emph{kernel} of a sequence 
$u=u_0u_1u_2\ldots$~\cite[Definition 2.2.18.]{Haeseler}, which (in one
dimension) one can think of as the set $\{u_au_{q^k+a}u_{2 q^k+a}\ldots
\mathbin| k\in\N_0, a\in\{0,\ldots,q^k-1\} \}$, i.e., the set of subsequence
on $q^k_{} \Z+a$ ($Q^k_{} L+a$) for all $k$ and $a$. The finiteness (and
cardinality) of this set is in question, and can be answered
(see~\cite[Lemma 2.2.22.]{Haeseler}) by the \emph{substitution
  graph}~\cite[Definition 2.2.20.]{Haeseler}: The substitution graph is the
directed, labelled graph with base point $(1,\ldots,m)$ (the whole
alphabet) and the vertices are given by all ordered $m$-tuples
$(\Phi_{}^k(1)_{z'},\ldots,\Phi^k_{}(m)^{}_{z'})$ (where $z'\in L/Q^k_{}L$),
which can be obtained by some power $k$ of the substitution. We draw 
an edge with label $z$ from 
$v^{}_1=(s^{}_1,\ldots,s^{}_m)$ to
$v^{}_2=(\Phi(s^{}_1)^{}_z,\ldots,\Phi(s^{}_m)^{}_z)$ 
($z\in L/QL$).

Obviously, in the case $L=L'$, the coincidence graph can be obtained
from the substitution graph: $m$-tuples which are equal as sets are
identified. And Lemma~\ref{lem:children} tells us that with this
identification there is no problem with the edges. Furthermore, an
$m-tuple$ $(i,\ldots,i)$ indicates a coincidence. As an example, we
give the substitution (see~\cite[Example 1. following Definition
2.2.20.]{Haeseler}) and coincidence graphs of the Thue-Morse sequence 
($a\mapsto ab$, $b\mapsto ba$).    

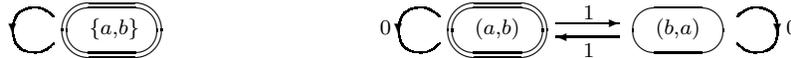
\begin{figure}[h]
\begin{tabular}{c@{\hspace*{2cm}}c}
\setlength{\unitlength}{.6cm}
\begin{picture}(5,2)
\qbezier(1.7,1.3)(1.5,1.5)(1.3,1.5)
\qbezier(.8,1)(.8,1.5)(1.3,1.5)
\qbezier(.8,1)(.8,0.5)(1.3,0.5)
\qbezier(1.7,0.7)(1.5,0.5)(1.3,0.5)
\put(.8,0.9){\vector(0,-1){0}}

\put(3,1){\oval(2,1)}
\put(3,1){\oval(2.2,1.2)}
\put(2.5,0.9){$\scriptstyle \{a,b\}$}
\end{picture}
&
\setlength{\unitlength}{.6cm}
\begin{picture}(10,2)
\put(.4,.9){$\scriptstyle 0$}
\qbezier(1.7,1.3)(1.5,1.5)(1.3,1.5)
\qbezier(.8,1)(.8,1.5)(1.3,1.5)
\qbezier(.8,1)(.8,0.5)(1.3,0.5)
\qbezier(1.7,0.7)(1.5,0.5)(1.3,0.5)
\put(.8,0.9){\vector(0,-1){0}}

\put(3,1){\oval(2,1)}
\put(3,1){\oval(2.2,1.2)}
\put(2.5,0.9){$\scriptstyle (a,b)$}

\put(4.3,1.15){\vector(1,0){1.4}}
\put(5.7,.85){\vector(-1,0){1.4}}
\put(4.9,1.25){$\scriptstyle 1$}
\put(4.9,.4){$\scriptstyle 1$}

\put(7,1){\oval(2,1)}
\put(6.5,0.9){$\scriptstyle (b,a)$}

\qbezier(8.3,1.3)(8.5,1.5)(8.7,1.5)
\qbezier(9.2,1)(9.2,1.5)(8.7,1.5)
\qbezier(9.2,1)(9.2,0.5)(8.7,0.5)
\qbezier(8.3,0.7)(8.5,0.5)(8.7,0.5)
\put(9.2,0.9){\vector(0,-1){0}}
\put(9.4,.9){$\scriptstyle 0$}
\end{picture}
\\
\end{tabular}
\caption{ \label{fig:thuem}
Coincidence graph (left) and Substitution graph (right)
of the Thue-Morse substitution}
\end{figure}

\end{appendix}

\section*{Acknowledgements}

Both authors thank the referees for many valuable suggestions. 
D.F.~thanks L.~Danzer for helpful discussions
about nonadmissible substitutions. 
B.S.~acknowledges financial support by the Cusanuswerk.
Both authors express their thanks
to the German Research Council, Collaborative Research Center 701.


\end{document}